\documentclass[onefignum,onetabnum]{siamonline250211}

\usepackage{lipsum}
\usepackage{amsfonts}
\usepackage{graphicx}
\usepackage{epstopdf}
\usepackage{algorithmic}
\ifpdf
  \DeclareGraphicsExtensions{.eps,.pdf,.png,.jpg}
\else
  \DeclareGraphicsExtensions{.eps}
\fi

\usepackage{enumitem}
\setlist[enumerate]{leftmargin=.5in}
\setlist[itemize]{leftmargin=.5in}

\newsiamremark{remark}{Remark}
\newsiamremark{hypothesis}{Hypothesis}
\crefname{hypothesis}{Hypothesis}{Hypotheses}
\newsiamthm{claim}{Claim}
\newsiamremark{fact}{Fact}
\crefname{fact}{Fact}{Facts}

\usepackage{nicefrac}
\usepackage{siunitx}

\usepackage{tikz}
\tikzset{every picture/.style={font issue={\fontsize{10}{10}}},
        font issue/.style={execute at begin picture={#1\selectfont}}
        }

\usepackage{pgfplots}
\pgfplotsset{compat=1.18, scaled ticks=false}
\usepackage{pgfplotstable}
\pgfplotsset{
    discard if/.style 2 args={
        x filter/.code={
            \edef\tempa{\thisrow{#1}}
            \edef\tempb{#2}
            \ifx\tempa\tempb
                
            \fi
        }
    },
    discard if not/.style 2 args={
        x filter/.code={
            \edef\tempa{\thisrow{#1}}
            \edef\tempb{#2}
            \ifx\tempa\tempb
            \else
                
            \fi
        }
    }
}
\pgfplotstableset{
    highlight bold/.style = {
        postproc cell content/.style={
            @cell content/.add={\boldmath}{},
        },
    },
}

\def\phitheta   {\varphi_{\theta}}

\newcommand{\N}{\mathbb{N}}

\newcommand{\R}{\mathbb{R}}

\headers{Invertible ResNets for Inverse Imaging Problems}{C. Arndt and J. Nickel}
\definecolor{carnelian}{rgb}{0.7, 0.11, 0.11}
\title{Invertible ResNets for Inverse Imaging Problems: Competitive Performance with Provable Regularization Properties\thanks{Equal contribution. Alphabetical author order.
\funding{This work was supported by the Deutsche Forschungsgemeinschaft (DFG) - project number 281474342/GRK2224/2 - as well as the Bundesministerium für Bildung und Forschung (BMBF) - funding code 03TR07W11A and 05M22LBA. The responsibility for the content of this publication lies with the authors.}}}

\author{Clemens Arndt\thanks{Center for Industrial Mathematics, University of Bremen, 28359 Bremen 
  (\email{carndt@uni-bremen.de, junickel@uni-bremen.de}).}
\and Judith Nickel\footnotemark[2]}

\usepackage{amsopn}

\ifpdf
\hypersetup{
  pdftitle={Invertible ResNets for Inverse Imaging Problems},
  pdfauthor={C. Arndt and J. Nickel}
}
\fi

\begin{document}

\maketitle
 
\begin{abstract}
Learning-based methods have demonstrated remarkable performance in solving inverse problems, particularly in image reconstruction tasks. Despite their success, these approaches often lack theoretical guarantees, which are crucial in sensitive applications such as medical imaging. Recent works by Arndt et al addressed this gap by analyzing a data-driven reconstruction method based on invertible residual networks (iResNets). They revealed that, under reasonable assumptions, this approach constitutes a convergent regularization scheme. However, the performance of the reconstruction method was only validated on academic toy problems and small-scale iResNet architectures. 
In this work, we address this gap by evaluating the performance of iResNets on two real-world imaging tasks: a linear blurring operator and a nonlinear diffusion operator.
To do so, we compare the performance of iResNets against state-of-the-art neural networks, revealing their competitiveness at the expense of longer training times. Moreover, we numerically demonstrate the advantages of the iResNet's inherent stability and invertibility by showcasing increased robustness across various scenarios as well as interpretability of the learned operator, thereby reducing the black-box nature of the reconstruction scheme.
\end{abstract}



\section{Introduction}
\label{sec:introduction}

Inverse problems arise in a wide range of applications, such as image processing, medical imaging, and non-destructive testing. These problems share a common objective: recovering an unknown cause from possibly noisy measurement data by inverting the measurement process (i.e., the forward operator). A major challenge is that the causes often depend discontinuously on the observed data, a characteristic that makes these problems ill-posed.
To address this, the problems need to be regularized in order to allow for stable reconstructions. However, this adjustment of the problem must be small to ensure that the solutions remain accurate.

The theory behind classical regularization schemes is well established, and provides strong guarantees regarding their regularization properties \cite{Burger_2018}. Recently,  deep learning based approaches have demonstrated a remarkable performance in solving inverse problems. Many of these methods are built on classical algorithms, e.g., unrolled iterative schemes \cite{monga2021} or plug-and-play methods \cite{Venkata2013}. However, the theoretical understanding in terms of provable guarantees remains limited. 
A key concern is the lack of stability in most approaches, which can pose significant risks in safety-critical applications. For this reason, there has been a growing body of research aimed at investigating regularization properties for deep learning methods (cf.\ ``related literature'' below).

We tackle the challenges of ill-posedness with a supervised learning approach based on invertible residual networks (iResNets). The basic idea is to approximate the forward operator with the forward pass of the network and use the inverse to solve the inverse problem. A general convergence theory, providing theoretical guarantees for iResNets, has been developed in \cite{iresnet_01_regtheory}. In addition, \cite{iresnet_02_bayesian} analyzed different training strategies and investigated to what extent iResNets learn a regularization from the training data. 
However, this approach has so far been limited to simplified toy examples, and we observed that the architecture proposed in \cite{iresnet_01_regtheory,iresnet_02_bayesian} struggles to scale effectively to real-world tasks. This raises a critical question: can iResNets be designed to achieve competitive performance with state-of-the-art methods on real-world applications?

To address this, we evaluate the performance and robustness of iResNets through extensive numerical experiments on both a linear blurring and a nonlinear diffusion problem, comparing the results with state-of-the-art methods.
The results reveal that iResNets can indeed achieve competitive performance, albeit with the trade-off of longer training times. However, this limitation can be partially compensated by the additional benefits resulting from the inherent stability and invertibility of the iResNet. Our numerical investigations show that these properties increase the robustness of the reconstruction method against adversarial attacks, varying noise levels during test time and limited training data. Moreover, the invertibility allows for an interpretation of the learned forward operator and its regularization, thereby reducing the ``black-box" nature commonly associated with learned regularization techniques.

This work thus provides a fully learned reconstruction scheme with strong regularization properties and high performance. 
Although several studies have investigated learning-based methods with theoretical guarantees -- such as approaches that learn regularization terms or employ plug-and-play techniques (cf.\ \cite{Mukherjee_2023}) -- there has been limited focus, to the best of our knowledge, on fully learned reconstruction methods with inherent regularization properties. Consequently, our work, in combination with \cite{iresnet_01_regtheory, iresnet_02_bayesian}, offers foundational contributions to this area and has the potential to inspire further research.

The manuscript is organized as follows. We begin in Section~\ref{sec:theory} by introducing iResNets and reviewing their regularization theory. Following this, in Section~\ref{sec:training_and_pde}, we elaborate why iResNets are well-suited for diffusion and blurring problems. In Section~\ref{sec:results}, we report our numerical results. First, we introduce the specific architecture design, offering valuable insights to guide the design of effective and efficient invertible architectures in Section~\ref{sec:architecture}.
Then, we compare the reconstruction quality of iResNets with several other common deep learning approaches and a classical reconstruction method in Section~\ref{sec:performance_comparison}. This is followed by a series of robustness studies, including adversarial attacks, variations in testing data, and limited training data availability in Section~\ref{sec:robustness}. 
Finally, we leverage the invertibility of the iResNet to compare the forward pass of the trained iResNet with the true forward operator and investigate the learned regularization in Section~\ref{sec:investigate_reg_prop}.
These analyses fully rely on the invertibility of iResNets, offering insights that are not attainable with other reconstruction methods. We conclude in Section~\ref{sec:outlook} with a summary of our findings and suggestions for future research directions.

\subsection*{Related literature}

In recent years, as neural networks have achieved remarkable performance, the focus on interpretability and providing theoretical guarantees for deep learning algorithms has gained increasing importance. This is reflected in the growing body of research dedicated to addressing these concerns.

A comprehensive overview of deep learning methods with regularization properties can be found in \cite{Mukherjee_2023}. Furthermore, \cite{obmann2023convergence} presents stability and convergence results for deep equilibrium models. To achieve these results, the learned component of the equilibrium equation must meet specific conditions, which can be satisfied using iResNets, for example. Based on the classical theory of variational regularization, \cite{dittmer2021, mukherjee_2021} introduced the adversarial convex regularizer. This approach relies on architecture restrictions which enforce strong convexity in order to obtain a stable and convergent regularization scheme.

One notable contribution aiming at enhancing the interpretability of the learned solution is the DiffNet architecture \cite{arridge_hauptmann_2020}, which is specifically designed to model potentially nonlinear diffusion processes. Inspired by the structure of numerical solvers, the network consists of multiple blocks representing discrete steps in time. Within each block, a differential operator, computed by a small subnetwork, is applied to the input. This design results in a highly parameter-efficient architecture that is fast to train and facilitates the analysis and interpretation of the differential operations.

In addition to research focused on providing theoretical guarantees or improving the interpretability of learned solutions, several approaches exhibit structural similarities to iResNets. One such method is the proximal residual flow \cite{hertrich2023}, a ResNet designed to be invertible, where the residual functions are constrained to be averaged operators, differing from the Lipschitz constraint we use in this work. This architecture has been used as a normalizing flow for Bayesian inverse problems. 
In \cite{ruthotto_2020}, convolutional ResNets are interpreted as discretized partial differential equations, deriving stability results for the architecture using monotonicity or Lipschitz conditions.

ResNet architectures are also employed in the context of Plug-and-play (PnP) algorithms.
In PnP frameworks, a denoising algorithm is incorporated to replace components of iterative algorithms used for solving inverse problems. However, this substitution typically leads to a loss of convergence guarantees, prompting extensive research to establish theoretical foundations for PnP algorithms.
Several studies leverage ResNet architectures with Lipschitz constraints -- similar to the proposed iResNet -- and have demonstrated both convergence guarantees and competitive performance (see \cite{sherry2024,ryu2019PnP,Hertrich2023PnP, Hurault2022PnP, Hurault2022PnP2}). These findings suggest that ResNet-style architectures are particularly well-suited for such applications. To ensure some degree of convergence in PnP algorithms, ResNets are often designed to be averaged operators. While the proposed iResNet does not inherently meet this requirement, it could be adapted following the approach in \cite[Lemma 18]{Hertrich2023PnP} by introducing an oracle to transform it into an averaged operator (see Appendix~\ref{sec:iResNet_averaged}).
Furthermore, recent work by \cite{Hurault2022PnP} has established convergence guarantees for a ResNet-style denoiser with a Lipschitz constraint on the residuals, similar to those in our iResNet framework, without necessitating the network to be averaged. 
Overall, these findings highlight the potential of iResNets as viable candidates for denoisers within PnP schemes, paving the way for further exploration in this area.

From a broader perspective, neural networks with stability guarantees have recently attracted significant attention. Particularly, constraining the Lipschitz constant of neural networks has become a common approach to enhance robustness, especially against adversarial attacks \cite{Gouk2021LipschitzReg, miyato2018spectral, Tsuzuku2018LipschitzAdversarial}.

We note that approaches which achieve both a high numerical performance and strong theoretical guarantees, such as stability and convergence, are still quite rare in general. Therefore, our iResNet approach provides an advantageous combination of highly desirable properties.

\section{iResNets for inverse problems}
\label{sec:theory}

We study iResNets $\varphi_\theta:X\to X$ with network parameters $\theta$ acting on a Hilbert space $X$, where 
\begin{equation}\label{Eq:Def concatenated iResNet}
    \varphi_\theta =  \varphi_{\theta_1,1} \circ ... \circ \varphi_{\theta_N,N}
\end{equation}
with $N\in \N$ and $\varphi_{\theta_i,i}:X\to X$ for $i\in \{ 1, \dots, N\}$. Each subnetwork $\varphi_{\theta_i,i}$ is defined as
\begin{equation*}
\varphi_{\theta_i,i} = \text{Id} - f_{\theta_i,i} \quad \text{for } i \in \{1, \dots, N\}
\end{equation*}
with Lipschitz continuous residual functions $f_{\theta_i,i}: X \to X$ satisfying $\text{Lip}(f_{\theta_i,i}) \leq L_i < 1$. The bound on the Lipschitz constant of $f_{\theta_i,i}$ allows for an inversion of the subnetworks $\varphi_{\theta_i,i}$ via Banach's fixed-point theorem with fixed-point iteration
\begin{equation*}
    x^{k+1} = z + f_{\theta_i,i}(x^k)
\end{equation*}
converging to $x = \varphi_{\theta_i,i}^{-1}(z)$ for $z\in X$. Moreover, each subnetwork $\varphi_{\theta_i,i}$ and its inverse $\varphi_{\theta_i,i}^{-1}$ satisfy the Lipschitz bounds
\begin{equation*}
    \text{Lip}(\varphi_{\theta_i,i}) \leq 1 + L_i \quad \text{and} \quad \text{Lip}(\varphi_{\theta_i,i}^{-1}) \leq \frac{1}{1-L_i}.
\end{equation*}
We refer the reader to \cite{iresnet_01_regtheory, behrmann2019invertible} for a detailed derivation of these results. The properties of the subnetworks directly translate to the concatenated network $\varphi_\theta$, making it invertible as each subnetwork is invertible with
\begin{equation}\label{Eq:Lip network and inverse}
    \text{Lip}(\varphi_{\theta}) \leq \prod_{i=1}^N (1 + L_i) \quad \text{and} \quad \text{Lip}(\varphi_{\theta}^{-1}) \leq \prod_{i=1}^N \frac{1}{1-L_i}.
\end{equation}
The Lipschitz constant of the inverse $\varphi_{\theta}^{-1}$ is crucial for ensuring the stability of the reconstruction scheme, as we will elaborate on later. To simplify our discussion, we define a Lipschitz parameter $L$ for the inverse $\varphi_{\theta}^{-1}$ such that
\begin{equation}\label{Eq:Lip concat inverse}
\text{Lip}(\varphi_{\theta}^{-1}) \leq \frac{1}{1-L} = \frac{1}{\prod_{i=1}^N(1-L_i)}.
\end{equation}

The aim is to apply the above defined iResNets to solve ill-posed inverse problems. Given that the input and output spaces of iResNets are inherently identical, we restrict our consideration to problems of the form
\begin{equation*}
    F(x) = z 
\end{equation*}
where $F \colon X \to X$ is a (possibly nonlinear) operator, and both the ground truth data $x\in X$ and the observed data $z\in X$ belong to the same Hilbert space $X$. 

Given the residual structure of iResNets, they are particularly well-suited for problems where learning deviations from the identity function is advantageous. This makes them ideal for applications such as image processing (e.g., denoising, deblurring), problems involving partial differential equations (PDEs), and pre-processing or post-processing tasks.

Furthermore, in certain scenarios, an operator $A$ mapping between different Hilbert spaces $X$ and $Y$ can be reformulated to fit the previously described setting by introducing a related operator $F$. More precisely, if the forward operator is linear -- that is, if we consider $A\in\mathcal{L}(X,Y)$ -- one can leverage the normal equation to transform the problem into a form that satisfies the previously stated conditions. Specifically, we define $F = A^\ast A$ with measurement data $z = A^\ast y = A^\ast A x$ for $x\in X$. For further details, we refer the reader to \cite[Remark 2.1]{iresnet_02_bayesian}.

Our objective is to recover the unknown ground truth $x^\dagger \in X$ via iResNets, given noisy measurements $z^\delta\in X$ with noise level $\delta > 0$ satisfying $\|z^\delta - F( x^\dagger)\| \leq \delta$. To do so, we follow the subsequent \textit{reconstruction approach}: \medskip
\begin{itemize}
    \item[(I)] The iResNet $\varphi_\theta:X\to X$ is trained to approximate the forward operator $F$. 
    \item[(II)] The inverse $\varphi_\theta^{-1}$ is used to reconstruct the ground truth $x^\dagger$ from $z^\delta$.
\end{itemize} \medskip

This reconstruction approach effectively addresses the ill-posedness of the inverse problem by satisfying essential properties of a regularization scheme. A regularization scheme aims to mitigate ill-posedness by employing a reconstruction algorithm that guarantees the existence and uniqueness of solutions, stability with respect to $z^\delta$ and convergence to the ground truth $x^\dagger$ as the noise level $\delta$ converges to zero.
The iResNet reconstruction approach is inherently designed to satisfy the first three conditions of a regularization scheme, namely existence, uniqueness, and stability of the reconstructed solutions.
This is mainly attributed to the Lipschitz continuity of the iResNet’s inverse.
To be more precise, the existence and uniqueness of $\varphi_\theta^{-1}(z^\delta)$ is synonymous with the invertibility of $\varphi_\theta$ and thus automatically fulfilled. Moreover, the stability, i.e., continuity of $\varphi_\theta^{-1}$, directly follows from Equation~\ref{Eq:Lip network and inverse}, cf.\ \cite[Lemma 3.1]{iresnet_01_regtheory}.

The only aspect not directly guaranteed is the convergence of $\varphi_{\theta}^{-1}(z^\delta)$ to the ground truth as $\delta\to 0$, which requires certain additional prerequisites.
Simply put, the convergence depends on the success of the training and the expressivity of the trained network. Notably, no additional assumptions about the forward operator or its null space are required. 
In \cite[Theorem 3.1]{iresnet_01_regtheory}, a convergence result has been derived in the setting of a linear forward operator and an iResNet consisting of a single subnetwork, i.e.\ $N=1$. 
However, this result can be directly extended to encompass nonlinear forward operators and iResNets with multiple subnetworks, as defined in Equation~\eqref{Eq:Def concatenated iResNet}.
For completeness, we present a simplified version of the convergence result from \cite[Theorem 3.1]{iresnet_01_regtheory}.
To do so, the network parameters $\theta$ must depend on the Lipschitz parameter $L$, as introduced in Equation~\ref{Eq:Lip concat inverse}. We emphasize this dependency by adapting the notation to $\varphi_{\theta(L)}$.

\begin{lemma}[Convergence, cf.\ \cite{iresnet_01_regtheory}]\label{lem:convergence}
For $x^\dagger\in X$, let $z^\delta \in X$ satisfy $\|z^\delta - F(x^\dagger)\| \leq \delta$. Moreover, assume that $x^\dagger$ and the network $\varphi_{\theta(L)}$ with network parameters $\theta(L)$ for $L \in [0,1)$ satisfy \begin{equation}\label{eq:inversion_error_convergence_prop}
    \| \varphi_{\theta(L)}^{-1}(F(x^\dagger)) - x^\dagger \|  \rightarrow 0 \quad \text{as } L \to 1.
\end{equation}
If the Lipschitz parameter $L:(0,\infty) \rightarrow [0,1)$ is chosen such that 
\begin{equation}\label{eq:conv_param_choice_ass}
    L(\delta) \rightarrow 1 \quad \text{and} \quad \frac{\delta}{1-L(\delta)} \rightarrow 0 \qquad \text{for } \delta \to 0,
\end{equation}
then it holds 
\begin{equation*}
    \| \varphi_{\theta(L(\delta))}^{-1}(z^\delta)  - x^\dagger\| \rightarrow 0 \qquad \text{for } \delta \to 0.
\end{equation*}
\end{lemma}

The proof can be found in Appendix~\ref{Sec:Proof_convergence}.

\begin{remark}
The convergence result established in \cite[Theorem 3.1]{iresnet_01_regtheory} relies on the so-called local approximation property
\begin{equation*}
   \|F(x^\dagger) -\varphi_{\theta(L)}(x^\dagger) \| = \mathcal{O}((1-L)\psi(1-L)) \quad \text{as } L \to 1,
\end{equation*}
where $\psi \colon \R_{\geq  0} \to \R_{\geq  0}$ is an index function, i.e.\ $\psi$ is continuous, strictly increasing and satisfies $\psi(0)=0$. In our setting, the local approximation property is replaced by Equation~\eqref{eq:inversion_error_convergence_prop}, which better aligns with the reconstruction training framework (cf.\ Section~\ref{sec:results}).
Note that both the local approximation property and Equation~\eqref{eq:inversion_error_convergence_prop} depend on the approximation capacity of the trained network at the specific element $x^\dagger$. Hence, unlike classical convergence results in the literature, the convergence in Lemma~\ref{lem:convergence} is local in nature. This characteristic bears similarities to classical source conditions. Specifically, the convergence result does not hold for all $x^\dagger\in X$, but only for those which the network can approximate with sufficient accuracy.
For a more detailed discussion and analysis of the local approximation property, we refer the reader to \cite{iresnet_01_regtheory}.
\end{remark}

\section{Appoximating diffusion processes with iResNets}
\label{sec:training_and_pde}

For our numerical experiments, we consider a nonlinear diffusion problem and a linear blurring operator (i.e., linear diffusion). These problem settings are particularly interesting, because several works observed an analogy between ResNets and numerical solvers for ODEs \cite{behrmann2019invertible, ruthotto_2020, sherry2024}. In the following, we briefly review this relationship by focusing on the invertibility of the iResNet-subnetworks.

Let $X = L^2(\Omega)$ be the space of ground truth images and data with some domain $\Omega \subset \R^n$ and $F \colon X \to X$ be an operator which maps a clean image to a diffused version of it. More precisly, we consider the forward problem $F(u_0) = u(T, \cdot)$ for some $T>0$, where $u \colon \R_{\geq 0} \times \Omega \to \R$ is the solution of the PDE
\begin{align}
    \partial_t u &= \mathrm{div} \left( g(| \nabla u |) \nabla u \right) & &\text{on } (0,T] \times \Omega,\label{eq:pm_diffusion} \\
    u(0, \cdot) &= u_0 & &\text{on } \Omega, \notag
\end{align}
known as Perona-Malik diffusion \cite{peronamalik1990}, with $g \in C^1(\R_{\geq 0})$ and zero Dirichlet or zero Neumann boundary condition on $\partial \Omega$. 

By considering the explicit Euler method
\begin{equation*}
    u_{t+1} = u_t + h \cdot \mathrm{div} \left( g(| \nabla u_t |) \nabla u_t \right)
\end{equation*}
for solving the PDE numerically, a similarity to the structure of a ResNet (identity plus differential operator) can be observed.
However, since the differential operator is not continuous with respect to $u_t \in L^2 (\Omega)$, fitting the contractive residual function $f_{\theta, i}$ of an iResNet-subnetwork to it might be challenging.

This changes when we consider the implicit Euler method 
\begin{equation}
\label{eq:pm_implicit_euler}
    u_{t+1} = u_t + h \cdot \mathrm{div} \left( g(| \nabla u_{t+1} |) \nabla u_{t+1} \right),
\end{equation}
which is an elliptic PDE for $u_{t+1}$. 
In the following, we demonstrate that under reasonable assumptions, the solution operator $S_{g,h} \colon L^2(\Omega) \to L^2(\Omega)$, $S_{g,h} (u_{t}) = u_{t+1}$ can be approximated by iResNet-subnetworks. 
For this purpose, we consider a Perona-Malik filter function $g_\lambda$ which fulfills
\begin{equation} \label{eq:pm_filter_monotone}
g_\lambda(\tilde{s})\, \tilde{s} \leq g_\lambda(s)\, s \qquad \text{for} \qquad 0 \leq \tilde{s} \leq s \leq \lambda.
\end{equation}
This holds true for the most common choices of $g_\lambda$ \cite{peronamalik1990}. Additionally, we have to restrict to images $u_t$, where $|\nabla u_t| \leq \lambda$ holds in most parts of the image (i.e., the majority of the image consists of smooth areas and weak edges). At first, we focus on images that do not contain any strong edges.

\begin{lemma}
\label{lem:Sh_fne}
For some $\lambda > 0$, let $g_\lambda \in C^1(\R_{\geq 0})$ be a Perona-Malik filter function, which fulfills \eqref{eq:pm_filter_monotone}. Let $u_t, v_t \in H^{1,2}(\Omega)$ satisfy $|\nabla u_t|, |\nabla v_t| \leq \lambda$. Then, it holds
\begin{equation}
\label{eq:pm_operator_lipschitz}
\|(\mathrm{Id} - S_{g_\lambda,h}) (u_t) - (\mathrm{Id} - S_{g_\lambda,h}) (v_t)\|^2 + \|S_{g_\lambda,h} (u_t) - S_{g_\lambda,h} (v_t)\|^2 \leq \|u_t - v_t\|^2.
\end{equation}
\end{lemma}

The proof can be found in Appendix~\ref{sec:proof_Sh_fne}.\\
This lemma shows that the mapping $\mathrm{Id} - S_{g_\lambda, h}$ has a Lipschitz constant of at most one on the set of images without strong edges. This implies that $S_{g_\lambda, h}$ can be accurately modeled with a subnetwork $\varphi_{\theta, i}$, as this requires $f_{\theta, i} \approx \mathrm{Id} - S_{g_\lambda, h}$.
However, if two images $u_t$, $v_t$ do contain strong edges, \eqref{eq:pm_operator_lipschitz} might be violated. But since natural images mainly consist of smooth regions and weak edges, with only a small proportion of strong edges, it is highly probable that at least
\begin{equation}
\|(\mathrm{Id} - S_{g_\lambda,h}) (u_t) - (\mathrm{Id} - S_{g_\lambda,h}) (v_t)\|^2 \leq \|u_t - v_t\|^2,
\end{equation}
i.e., $\mathrm{Lip}(\mathrm{Id} - S_{g_\lambda,h}) \leq 1$ is fulfilled. Therefore, we can expect $\varphi_{\theta, i}$ to be able to approximate $S_{g_\lambda, h}$ on the manifold of natural images. 
We refrain from delving into the specific number of strong edges permissible, as this depends on the parameters $g_\lambda$ and $\lambda$. Instead, this result may serve as a motivation for approximating diffusion processes for natural images with iResNets.

\section{Numerical Results}
\label{sec:results}

To validate the performance of iResNets in solving inverse problems, we conduct numerical experiments on two forward operators. As discussed in Section~\ref{sec:training_and_pde}, diffusion and blurring operators are particularly suitable for this purpose. Therefore, we consider a linear blurring operator defined as a convolution with a Gaussian kernel with a kernel size of $11\times 11$ and standard deviation of
$\nicefrac{5}{3}$. Additionally, we examine an anisotropic nonlinear diffusion problem governed by the PDE in Equation~\eqref{eq:pm_diffusion}, utilizing the Perona-Malik filter function $g(| \nabla u |) = \nicefrac{1}{1+\lambda^{-2}|\nabla u|^2}$ with contrast parameter $\lambda=0.1$. The diffused image is obtained after 5 steps of Heun's method with step size $0.15$ and zero Neumann boundary conditions. 

All networks are trained on pairs of distorted and clean grayscale images from the STL-10 dataset \cite{coates2011}, where we use $\num{16384}$ images for training, $\num{64}$ for validation and $\num{128}$ for testing and evaluation. The distorted images are generated by applying the forward operator to the clean images, combined with additive Gaussian white noise with standard deviation $\delta>0$. We consider three different levels of noise in our numerical experiments, namely $\delta=0.01$, $\delta=0.025$, and $\delta=0.05$.

In \cite{iresnet_02_bayesian}, two fundamentally different approaches for training an iResNet to solve an inverse problem are discussed. Both methods require supervised training with paired data $(x_i, z_i^\delta)$. The first approach, referred to as \textit{approximation training}, aims at training $\varphi_\theta$ to approximate $F$ via the training objective
\begin{equation*}
    \min_{\theta} \sum_i \| \varphi_\theta(x_i) - z_i^{\delta} \|^2.
\end{equation*}
The second approach, referred to as \textit{reconstruction training}, focuses on training the inverse $\varphi_\theta^{-1}$ to reconstruct the ground truth via the training objective
\begin{equation} \label{eq:reconstruction_training}
    \min_{\theta} \sum_i \| \varphi_\theta^{-1}(z_i^{\delta}) - x_i \|^2.
\end{equation}
Reconstruction training is preferable for obtaining a data-driven regularization method, since it results in an approximation of the posterior mean estimator of the training data distribution \cite[Lemma 4.2]{iresnet_02_bayesian}. In contrast, approximation training overlooks many features of the training data distribution, with its regularizing effect primarily arising from the architectural constraints \cite[Theorem 3.1]{iresnet_02_bayesian}.  For this reason, we perform reconstruction training in our numerical experiments. It is important to note that this approach requires computing the network's inverse during training, which is done via a fixed-point iteration within each subnetwork, cf.\ Section~\ref{sec:theory}. 

The source code for the experiments in this section is available on GitLab\footnote{\url{https://gitlab.informatik.uni-bremen.de/junickel/iresnets4inverseproblems}}. 

\subsection{Architecture}
\label{sec:architecture}

\begin{figure}[t]

\centering
\includegraphics[width=0.99\textwidth]{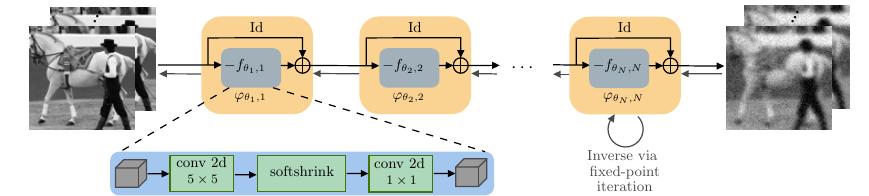}
\caption{Architecture of the iResNets used in our numerical experiments.}\label{fig:architecture}
\end{figure}

In this section, we discuss key practical considerations in the design of iResNets to achieve high-performing networks, along with details of the architectures used in our numerical experiments.
We tested a small and a big architecture (denoted as ``small iResNet'' and ``iResNet'', respectively), both designed according to \eqref{Eq:Def concatenated iResNet}.
Each residual function $f_{\theta_i, i}$ is defined as a 5$\times$5-convolution followed by a soft shrinkage activation function with learnable threshold and a 1$\times$1-convolution. While most common activation functions (e.g., ReLU, leaky ReLU, sigmoid) have a Lipschitz constant of one, there often exist large regions, where the actual slope is significantly smaller (e.g., for negative inputs of ReLU). As we are interested in estimating the Lipschitz constant of each part of the network as tightly as possible, soft shrinkage activation is a good choice because it attains its maximal slope of one for all inputs whose absolute value is larger than the threshold of the function. More details about both architectures are listed in Table \ref{tab:network_params} and a visualization of the network design is depicted in Figure~\ref{fig:architecture}.

We made five different choices for the Lipschitz parameter of the architectures, i.e., $L=0.95$, $L=0.97$, $L=0.99$, $L=0.995$, and $L=0.999$, to be able to test different levels of stability. Unless otherwise stated, the reported results correspond to the network with the highest $L$ because it is the least restricted. To control the Lipschitz constants of the residual functions $f_{\theta_i, i}$, we need to regulate the operator norms of each convolution, which can be efficiently computed with power iterations.
Thus, before performing a convolution, we apply a projection to the (trainable) network parameters to obtain convolutional weights which fulfill the desired operator norm.
It is important to note that during training, automatic differentiation accounts for both the projection step and the computation of operator norms. This strategy was already used in \cite{iresnet_02_bayesian} and it prevents the gradient steps computed by the optimizer from being in conflict with the Lipschitz constraint.

Previous studies \cite{iresnet_01_regtheory, iresnet_02_bayesian} introduced and evaluated an iResNet architecture on small-scale toy problems. In these works, the proposed architecture consists of an iResNet with a single subnetwork ($N=1$). Additionally, the authors suggest constructing the inverse of the iResNet based on their observation that it can be directly implemented as a scaled iResNet \cite[Remark 2.3]{iresnet_02_bayesian}. This approach eliminates the need for fixed-point iterations during training, potentially reducing training time.

However, upon evaluating both of these architectural choices, we found that they did not scale well to real-world tasks. Instead, we adopted an alternative approach, incorporating a larger number of subnetworks in our architectures -- $N=12$ for the smaller model and $N=20$ for the larger -- while keeping the residual functions $f_{\theta,i}$ relatively shallow. The rationale behind this design choice is discussed in the following sections.

\begin{table}[t]
    \centering
        \caption{Network and training parameters of our iResNets trained on a system with an AMD EPYC 7742 64-Core processor and a GeForce RTX 3090 GPU.} 
    \label{tab:network_params}
        \begin{tabular}{ c | c c c c c}
          & iResNet & Small iResNet & DiffNet & U-Net & ConvResNet \\ \hline
         Subnetworks/Scales & 20 & 12 & 5 & 5 & 20\\  
         Input channels & 64 & 32 & 1 & 1 & 64\\
         Hidden channels & 128 & 64 & 32 & 16--256 & 128\\
         Kernel size & 5 & 5 & 3 & 5 & 5 \\
         Trainable parameters & $\num{4263700}$ & $\num{640140}$ & \num{101310} & \num{5468705} & \num{4263873} \\ \hline
         Learning rate & $\num{1e-4}$ & $\num{1e-3}$ & $\num{2e-3}$ & $\num{5e-4}$ & $\num{1e-4}$\\
         Epochs & 500 & 300 & 100 & 500 & 300 \\
         Batch size & 32 & 32 & 16 & 64 & 16\\
         Optimizer & Adam & Adam & Adam & Adam & Adam\\
         Training time & 6--7 days & 1--2 days & 1 hour & 2 hours & 15 hours
        \end{tabular}
\end{table}

First, the Lipschitz constant of a neural network can be estimated from above by the product of the Lipschitz constants of all layers, but this estimation can be far from tight, as noted in \cite{bungert_2021}. This issue gets worse as the number of layers increases, which is why we chose shallow sub-architectures to mitigate this problem.

Second, if $L_i$ is chosen to be very close to one, the invertibility condition $\mathrm{Lip}(f_{\theta, i}) < 1$ might in fact be violated due to small numerical errors. This can lead to severe problems when attempting to invert $\varphi_{\theta,i}$. 
By using a large number of subnetworks $N$, we can assign smaller, more stable values to $L_i$ while still achieving a large overall $L$ value (see \eqref{Eq:Lip concat inverse}). For our small architecture with $N=12$, we can, e.g., choose $L_i \approx 0.438$ to obtain $L=0.999$. 

Third, the residual functions $f_{\theta, i}$ must be evaluated in each fixed-point iteration for the inversion of $\varphi_{\theta, i}$. Smaller subnetworks make these computations faster. Additionally, smaller $L_i$ values reduce the number of iterations required to reach a specific accuracy level. Although a higher $N$ increases the number of subnetworks to invert, the overall inversion process is faster with shallow residual functions. 

Fourth, a single iResNet-subnetwork has a Lipschitz constant of at most $1+L_i < 2$ in the forward pass. 
This implies that $\|\varphi_{\theta, i}^{-1}(z_1) - \varphi_{\theta, i}^{-1}(z_2)\| \geq \frac{1}{2} \|z_1 - z_2\|$, which limits the regularization capabilities of $\varphi_{\theta, i}^{-1}$ since it cannot map different data $z_1, z_2$ arbitrarily close to the same solution $x^\dagger$. Increasing the number of subnetworks alleviates this issue.

Fifth, while individual subnetworks are always monotone (due to $\langle \varphi_{\theta, i} (x), x \rangle = \|x\|^2 - \langle f_{\theta, i} (x), x \rangle \geq (1-L_i) \|x\|^2$), concatenating several of them allows for fitting increasingly non-monotone functions, thereby expanding the set of suitable forward operators.

Finally, as we consider the solution operator of a PDE (see Section \ref{sec:training_and_pde}) as the forward operator for our numerical experiments, a concatenated architecture is a natural choice. 

A drawback of using multiple concatenated shallow subnetworks is that the number of channels (i.e., the network's width) can only be expanded within the residual functions. The input and output dimension of all subnetworks has to be the same. To overcome this limitation, we lift the inverse problem $F \colon X \to X$ to a multichannel problem $\tilde{F} \colon X^M \to X^M$, $\tilde{F}(x_j)_k = F(x_j)$ for $j,k=1, ..., M$. Input and target images $z^\delta, x^\dagger \in X$ for training and evaluation of $\varphi_\theta$ are accordingly stacked to multichannel representations $(z^\delta, ..., z^\delta), (x^\dagger, ..., x^\dagger) \in X^M$. The mapping of the network's output back to the original space $X$ is simply implemented as a mean over all channels. Note that in this setting $\varphi_\theta$ is not an invertible mapping on $X$ but on $X^M$. 

To perform reconstruction training \eqref{eq:reconstruction_training}, the gradients of $\varphi_\theta^{-1}$ are required. We use the idea of deep equilibrium models \cite{bai_kolter_2019} to avoid the memory-intensive process of backpropagation through the potentially large number of fixed-point iterations. 
This means that the derivative of the solution $x$ of
\begin{equation*}
    x - f_{\theta, i}(x) = z
\end{equation*}
with respect to $\theta$ and $z$ is directly calculated via the implicit function theorem. In \cite{Gilton2021DeepEA}, this technique has already been used for solving inverse problems to simulate an infinite number of layers in unrolled architectures.

\subsection{Performance and comparison to other models}
\label{sec:performance_comparison}

In this section, we assess the performance of the iResNet reconstruction scheme and compare it with both traditional and deep-learning-based reconstruction methods.

\begin{figure}[p]
\centering
\includegraphics[width=0.96\textwidth]{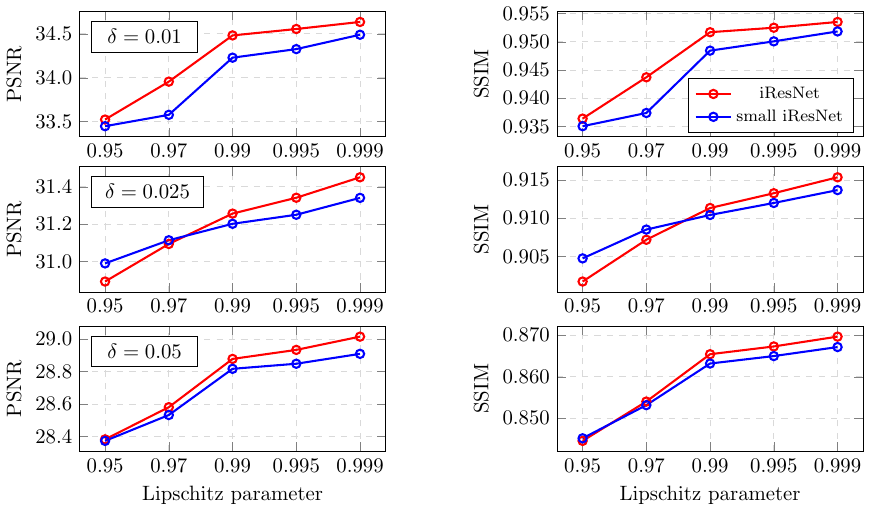}
\caption{Reconstruction performance of iResNets dependent on the Lipschitz parameter $L$ for the \textbf{nonlinear diffusion operator} and different noise levels ($\delta = 0.01, 0.025, 0.05$).}\label{fig:performance depending on L, nonlinear}

\vspace{8mm}

\includegraphics[width=0.96\textwidth]{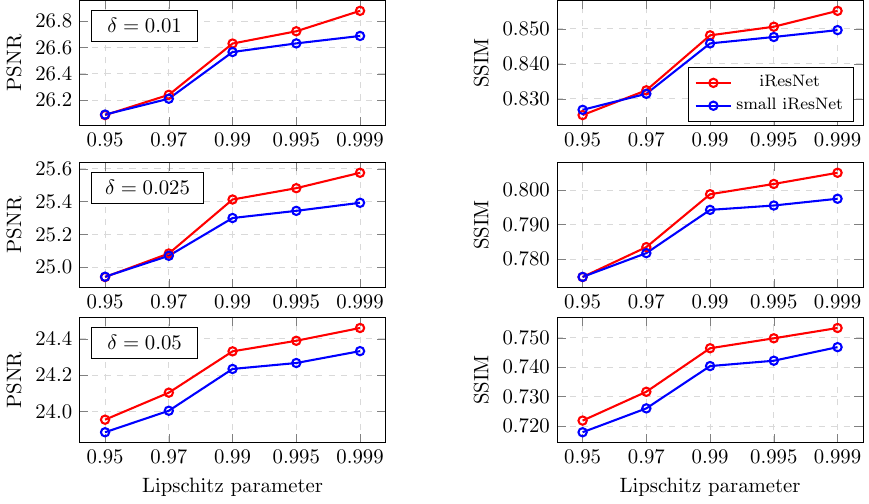}
\caption{Reconstruction performance of iResNets dependent on the Lipschitz parameter $L$ for the \textbf{linear blurring} and different noise levels ($\delta = 0.01, 0.025, 0.05$).}\label{fig:performance depending on L, linear}
\end{figure}

The performance of the iResNet and the small iResNet with respect to PSNR and SSIM across all three noise levels $\delta$ depending on the Lipschitz parameter $L$ is illustrated in Figure~\ref{fig:performance depending on L, nonlinear} for the nonlinear diffusion operator and in Figure~\ref{fig:performance depending on L, linear} for the linear blurring operator. One can observe that the reconstruction performance increases as the Lipschitz parameter $L$ approaches one for both forward operators and across all noise levels. Additionally, the larger iResNet outperforms the smaller iResNet at Lipschitz parameters close to one, while their performance is similar on average at lower Lipschitz parameters. This behavior is expected, as the Lipschitz constraint reduces the network's expressiveness, particularly at smaller Lipschitz parameters. Based on these observations, the following investigations will focus on the highest Lipschitz parameter, $L=0.999$. 

To further assess the reconstruction quality of our iResNets, we compare their performance against three other deep learning methods and the classical TV regularization. Specifically, we implement a convolutional ResNet (ConvResNet) without Lipschitz constraints on the residual functions, exhibiting the same network architecture as our iResNet, but with an additional learnable $1\times 1$ convolution at the beginning and end of the network. This modification overcomes the need to transfer the problem to a multichannel problem as discussed in Section~\ref{sec:architecture}.
Additionally, we employ the widely recognized U-Net architecture \cite{ronneberger2015} as well as the DiffNet from \cite{arridge_hauptmann_2020}, a specialized architecture based on non-stationary filters designed specifically for diffusion problems. We optimize the hyperparameters of our iResNets, TV, ConvResNet and U-Net on the validation set using random search. The hyperparameters of DiffNet are adapted from \cite{arridge_hauptmann_2020}, as they also focus on a linear blurring and nonlinear diffusion forward operator. The architecture, training parameters, and corresponding training times for all networks are detailed in Table~\ref{tab:network_params}. One can observe that the DiffNet has the fewest trainable parameters, followed by the small iResNet, while the (larger) iResNet, U-Net, and ConvResNet have a similar and significantly higher number of parameters. Moreover, there is a considerable difference in training times: DiffNet and U-Net converge within 1--2 hours, whereas our iResNets require several days up to a week to converge.
This is mainly caused by two aspects. First, controlling the Lipschitz constant of every Layer requires extra computations in every training step (cf.\ Section~\ref{sec:architecture}). Second, the computation of the inverse $\varphi_\theta^{-1}$ is done via fixed-point iterations, making the evaluation of $\varphi_\theta^{-1}$ significantly more expensive than the evaluation of $\varphi_\theta$. We would like to stress that a fixed number of iterations is used for both computing the Lipschitz constant and performing the inversion. 

The reconstruction quality in terms of PSNR and SSIM for all methods at different noise levels is reported in Table~\ref{tab:performance_all_models_diffusion} for the nonlinear diffusion and in Table~\ref{tab:performance_all_models_linblur} for the linear blurring. 
The networks substantially outperform the classical TV reconstruction, with the ConvResNet achieving the highest SSIM and PSNR values across both forward operators. Our iResNet exhibits a similar performance than DiffNet and U-Net. Notably, for the nonlinear diffusion operator, the iResNet performs particularly well at high noise levels, while for the linear blurring, it surpasses both DiffNet and U-Net at all noise levels. The small iResNet achieves slightly lower reconstruction quality compared to the larger iResNet, but it still significantly outperforms TV reconstruction. 

\begin{table}[t]
\centering
\caption{Comparison of the reconstruction quality (PSNR, SSIM) between the iResNet architectures (with $L=0.999$) and other models for the \textbf{nonlinear diffusion operator} across various noise levels.}
\label{tab:performance_all_models_diffusion}
    \begin{filecontents*}{csv_files/all_models_diffusion.csv}
model,psnr-01,ssim-01,psnr-025,ssim-025,psnr-05,ssim-05
ConvResNet,34.92778796354096,0.957236485398708,31.672218620898526,0.9213016594913472,29.17307996935976,0.8762490672189067
U-Net,34.658252550188735,0.9551626805869148,31.46832788482881,0.9181346427265376,29.017201013086908,0.8722828401733631
DiffNet,34.67434458677523,0.9545188981707404,31.474065135362967,0.9169631635956178,28.971289465093083,0.8696557122948407
iResNet,34.63643230925668,0.953523706593062,31.45108134216875,0.9153446503966488,29.015362892705543,0.8696729332824485
Small iResNet,34.490510251899984,0.951843256205809,31.3405622243386,0.9136767793250166,28.90941917211785,0.8671602437613292
TV,33.52535294018038,0.9430920913587266,30.265611600033285,0.8959350198011753,27.70289767661652,0.8352907164228289
\end{filecontents*}
\pgfplotstabletypeset[
    col sep=comma,
    columns/model/.style={string type, column name={ }, column type=l, 
    		},
    columns/psnr-01/.style={column name={PSNR}, precision=2, column type/.add={|}{}},
    columns/ssim-01/.style={column name={SSIM}, precision=3, column type/.add={}{|}},
    columns/psnr-025/.style={column name={PSNR}, precision=2},
    columns/ssim-025/.style={column name={SSIM}, precision=3, column type/.add={}{|}},
    columns/psnr-05/.style={column name={PSNR}, precision=2},
    columns/ssim-05/.style={column name={SSIM}, precision=3},
    assign column name/.style={/pgfplots/table/column name={\textbf{#1}}},
    every head row/.style={before row={ & \multicolumn{2}{c|}{$\delta = 0.01$} & \multicolumn{2}{c|}{$\delta = 0.025$} & \multicolumn{2}{c}{$\delta = 0.05$}  \\}, after row=\hline},
	every row no 2/.style={after row=\hline},
	every row no 4/.style={after row=\hline},
    every row 0 column 1/.style={highlight bold},
    every row 0 column 2/.style={highlight bold},
    every row 0 column 3/.style={highlight bold},
    every row 0 column 4/.style={highlight bold},
    every row 0 column 5/.style={highlight bold},
    every row 0 column 6/.style={highlight bold},
	fixed zerofill=true
    ]{csv_files/all_models_diffusion.csv}

\vspace{4mm}

\caption{Comparison of the reconstruction quality (PSNR, SSIM) between the iResNet architectures (with $L=0.999$) and other models for the \textbf{linear blurring operator} across various noise levels.}
\label{tab:performance_all_models_linblur}
    \begin{filecontents*}{csv_files/all_models_lin_blur.csv}
model,psnr-01,ssim-01,psnr-025,ssim-025,psnr-05,ssim-05
ConvResNet,27.064906438687167,0.8602731460021477,25.709709646346425,0.8123846251256119,24.567805337655315,0.7587262210590338
U-Net,26.67097580414596,0.8520743828509804,25.49349427278582,0.8022429705610293,24.406467899802344,0.7488329381746148
DiffNet,26.784552919463778,0.852326494552804,25.5070644367199,0.8032622618647552,24.40131153488268,0.750485599846026
iResNet,26.87920411180184,0.8551115277089519,25.575883400961096,0.8050933056955905,24.459643556658293,0.7533102813237942
Small iResNet,26.688087745230145,0.8496087764121082,25.392653669625336,0.7975357561183958,24.33224282717101,0.74680630845604
TV,25.81562030059105,0.8243770820567732,24.606388716389624,0.7676824543032995,23.588075773289408,0.7064013403070908
\end{filecontents*}
\pgfplotstabletypeset[
    col sep=comma,
    columns/model/.style={string type, column name={ }, column type=l, 
    		},
    columns/psnr-01/.style={column name={PSNR}, precision=2, column type/.add={|}{}},
    columns/ssim-01/.style={column name={SSIM}, precision=3, column type/.add={}{|}},
    columns/psnr-025/.style={column name={PSNR}, precision=2},
    columns/ssim-025/.style={column name={SSIM}, precision=3, column type/.add={}{|}},
    columns/psnr-05/.style={column name={PSNR}, precision=2},
    columns/ssim-05/.style={column name={SSIM}, precision=3},
    assign column name/.style={/pgfplots/table/column name={\textbf{#1}}},
    every head row/.style={before row={ & \multicolumn{2}{c|}{$\delta = 0.01$} & \multicolumn{2}{c|}{$\delta = 0.025$} & \multicolumn{2}{c}{$\delta = 0.05$}  \\}, after row=\hline},
	every row no 2/.style={after row=\hline},
	every row no 4/.style={after row=\hline},
    every row 0 column 1/.style={highlight bold},
    every row 0 column 2/.style={highlight bold},
    every row 0 column 3/.style={highlight bold},
    every row 0 column 4/.style={highlight bold},
    every row 0 column 5/.style={highlight bold},
    every row 0 column 6/.style={highlight bold},
	fixed zerofill=true
    ]{csv_files/all_models_lin_blur.csv}
\end{table}

\begin{figure}[p]
\centering
\includegraphics[width=\textwidth]{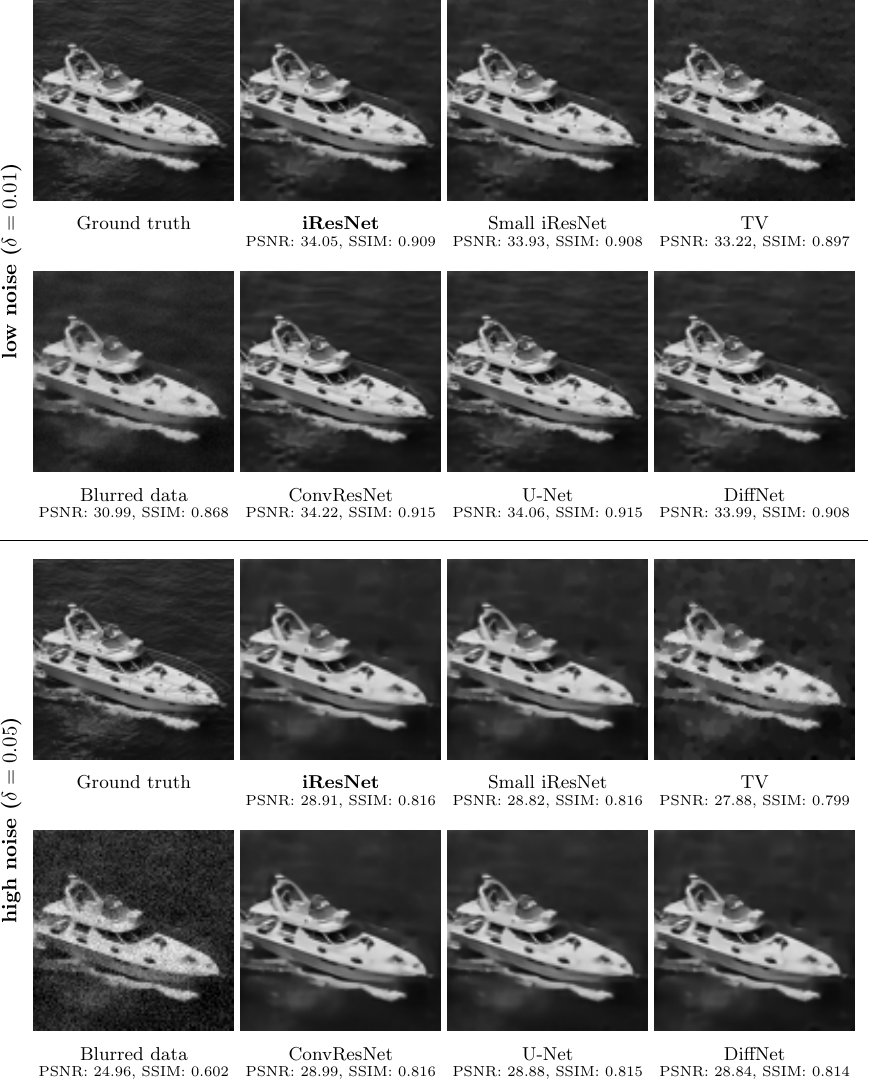}
\caption{Comparison of reconstructions from different methods for the \textbf{nonlinear diffusion operator}.}\label{fig:recos:nonlinear}
\end{figure}

\begin{figure}[p]
\centering
\includegraphics[width=\textwidth]{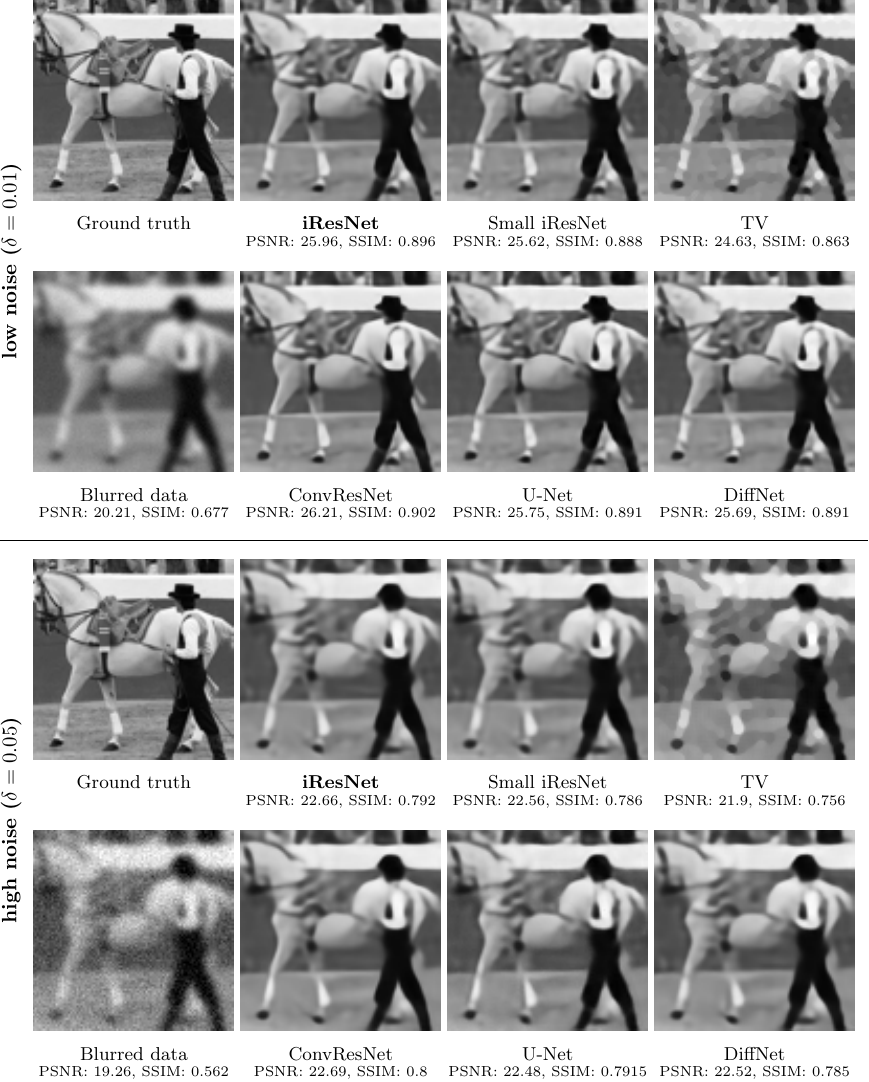}
\caption{Comparison of reconstructions from different methods for the \textbf{linear blurring operator}.}\label{fig:recos_linear}
\end{figure}

The quality of reconstructions produced by our iResNets is further validated through comparative examples, as can be seen in Figure~\ref{fig:recos:nonlinear} for the nonlinear diffusion and in Figure~\ref{fig:recos_linear} for the linear blurring at both the lowest and highest noise level. For both types of forward operators, the visual performance of our iResNet closely matches that of the ConvResNet, U-Net, and DiffNet. Moreover, the TV reconstructions exhibit severe artifacts, especially in the case of high noise and with the linear blurring operator. As anticipated, all methods struggle to recover fine details at high noise levels.

These results demonstrate that our iResNets perform on par with state-of-the-art methods in the literature, albeit with a significantly longer training time due to their invertibility. Among the compared methods, DiffNet appears to strike the best balance between the number of trainable parameters, training time, and performance. However, DiffNet is specifically designed for diffusion problems, whereas our iResNets are more versatile and can be applied to a broader range of forward operators, cf.\ Section~\ref{sec:training_and_pde}. 

We emphasize that iResNets are not expected to outperform all state-of-the-art methods, particularly since most competing approaches impose no explicit restrictions on network parameters, allowing for greater expressiveness. This raises a key question: Why choose iResNets if they only match the performance of state-of-the-art methods while requiring significantly longer training times? To address this, we present several studies demonstrating how the stability and invertibility of iResNets improve both the robustness and interpretability of the learned reconstruction process.

Given the superior performance of our larger iResNet (with $L=0.999$) compared to its smaller counterpart, we will focus our subsequent investigations on the larger model -- unless otherwise specified.

\subsection{Why stability matters}\label{sec:robustness}

The theoretical stability guarantee of iResNets also translates into practical advantages.
In comparison to other deep learning methods, iResNets exhibit greater robustness against adversarial data, higher noise levels during testing, and a significantly reduced number of training samples.
The following three subsections illustrate these benefits in detail.
Additionally, Appendix~\ref{sec:ood_data} provides an evaluation of the performance of all methods on out-of-distribution data.

\subsubsection{Adversarial examples}

\begin{figure}[t]
\includegraphics[width=0.95\textwidth]{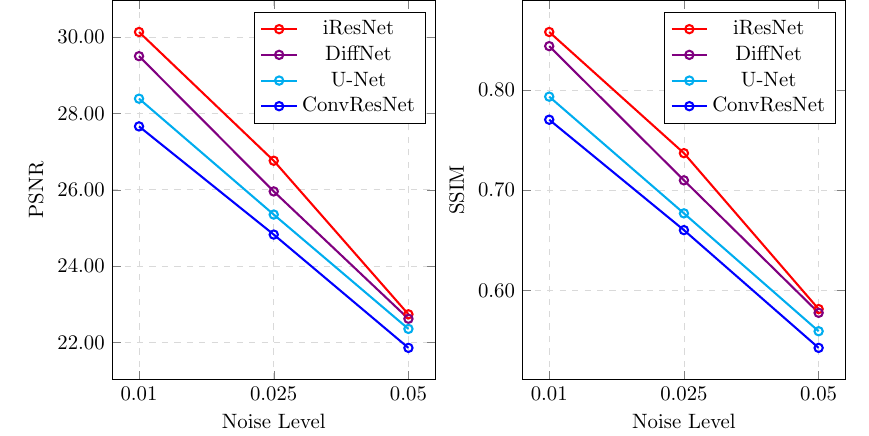}
\caption{Evaluation of the networks on adversarial examples for the \textbf{nonlinear diffusion operator}. For each of the four networks, we created five adversarial data samples by choosing the noise in a way that results in a particularly large reconstruction error. Then each network is evaluated on all 20 adversarial examples, and the average performance metrics are reported.}
\label{fig:eval_on_adv_examples}
\end{figure}


Adversarial data consists of distorted images that are created artificially for a certain reconstruction method $M$, in order to cause large reconstruction errors for $M$ (cf.\ \cite{szegedy2014}). More precisely, the measurement noise on the image $z^\delta$ is optimized via gradient ascent to maximize $\| M(z^\delta) - x^\dagger\|$ while keeping $\|z^\delta - F(x^\dagger)\| \leq \delta$. This way, $M$ can be tested in a worst-case scenario (regarding the noise).

We computed adversarial examples for the iResNet (with $L=0.999$), DiffNet, U-Net, and ConvResNet. Each model was then evaluated on the full set of adversarial examples -- including those generated using all other models, not just its own. 

For the nonlinear diffusion operator, the results in terms of PSNR and SSIM values of the reconstructions are illustrated in Figure~\ref{fig:eval_on_adv_examples}. A comparison with Table~\ref{tab:performance_all_models_diffusion} shows that the reconstruction error for all models is significantly higher on adversarial examples than on the normal test dataset.
The iResNet is the best-performing model on all noise levels, which can be attributed to its stability property. Specifically, its invertibility and stability ensure that the iResNet is well-defined across the entire measurement data space, enhancing its robustness against adversarial examples that lie outside the valid range of the forward operator.
For the linear blurring operator, one obtains similar results when choosing $L=0.995$, i.e., enforcing more stability.

\subsubsection{Test-time noise level variations}

\begin{figure}[p]
\centering
\includegraphics[width=0.99\textwidth]{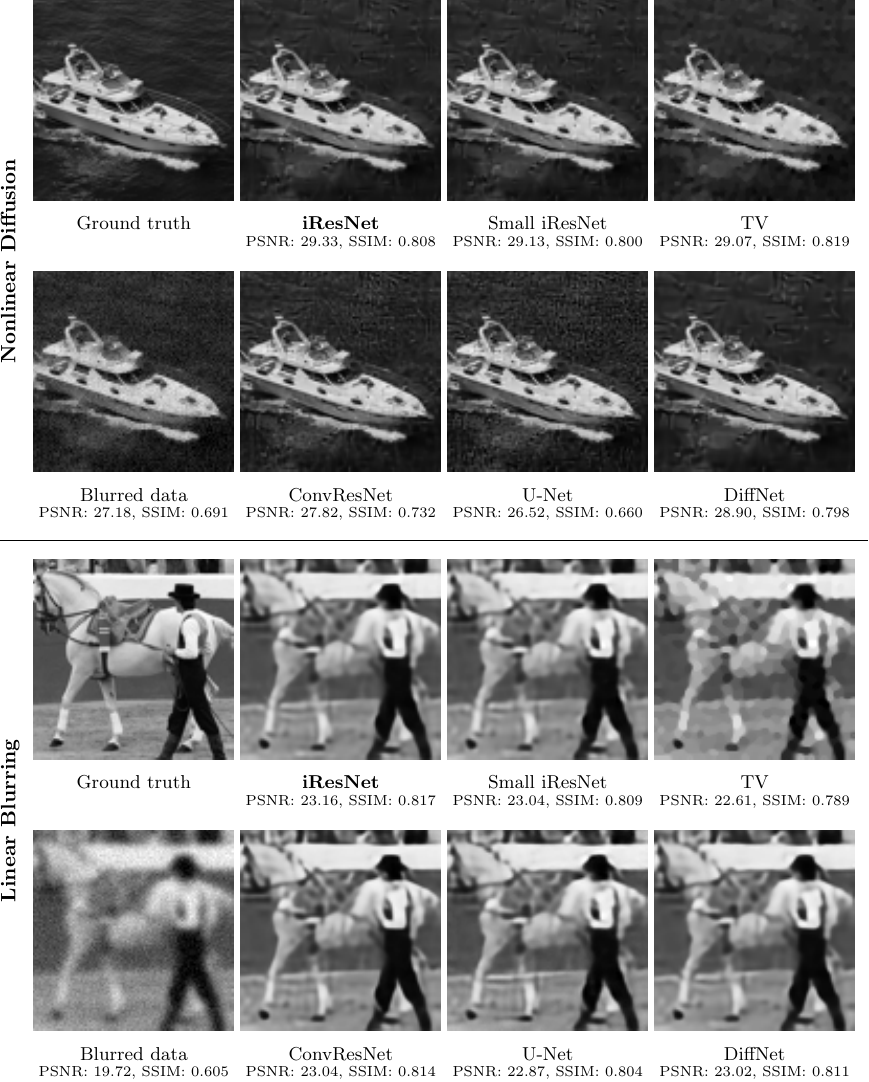}
\caption{Comparison of reconstructions for both forward operators using measurement data with noise level $\delta_{\text{test}}=1.4\cdot 0.025$.  All networks and hyperparameters were optimized on training data with $\delta=0.025$.}\label{fig:robustness_higher_noise}
\end{figure}

To analyze how the reconstruction methods respond to varying noise levels during testing, we consider the trained networks and TV, which were trained on the noise level $\delta = 0.025$, and evaluate them on testdata with noise level $\delta_\text{test}=c\, \delta$ with $c\in [0.7,1.8]$.
Figure~\ref{fig:robustness_higher_noise} presents reconstructions for both forward operators when the noise level is increased by a factor of $c=1.4$ compared to the training and hyperparameter tuning setup. 

The results indicate that the iResNet, DiffNet, and TV exhibit greater stability against noise level variations than the U-Net and ConvResNet. Notably, both iResNet variants outperform the other methods in terms of both visual quality and PSNR, with this effect being particularly pronounced for the nonlinear diffusion operator. In the nonlinear diffusion setting, ConvResNet and U-Net even appear to amplify noise in their reconstructions.

Additionally, we find that reducing the Lipschitz constant $L$ of the iResNet enhances robustness of the reconstructions against noise level changes, resulting in a behavior similar to TV for $L=0.95$. For an overview of the performance across the full range $c\in [0.7,1.8]$, we refer the reader to Appendix~\ref{App:varying_nlvl}.

\subsubsection{Small training dataset}\label{sec:small_dataset} 

We further assessed the performance of all models when trained on a significantly smaller dataset of 64 data pairs. In this setting, the risk of overfitting increases, making it more challenging for models to generalize to unseen test data.

Given its inherent stability and invertibility, the iResNet is expected to exhibit greater robustness under these conditions. The stability constraint serves as an implicit regularization mechanism during the training process, helping to mitigate overfitting. 

Table~\ref{tab:performance_fewdata64} presents the performance of all models across the nonlinear diffusion and linear blurring operators for the noise level $\delta = 0.025$. Notably, in this limited-data setting, using a Lipschitz parameter of $L=0.99$ or $L=0.995$ leads to better results compared to $L=0.999$, likely due to higher stability.

The iResNet performs on par with DiffNet, which is data efficient by construction \cite{arridge_hauptmann_2020}. In contrast, ConvResNet and U-Net exhibit a significantly weaker reconstruction performance, highlighting their tendency to overfit in low-data regimes. Additionally, we note that reducing the number of training samples also shortens the iResNet’s training time to 2–4 days.

\begin{table}[t]
\centering
    \caption{Reconstruction performance on the test data with noise level $\delta = 0.025$ using networks trained on 64 samples. The results on the left-hand side correspond to the \textbf{nonlinear diffusion operator} and on the right-hand side to the \textbf{linear blurring operator}.}
    \label{tab:performance_fewdata64}
    \begin{filecontents*}{csv_files/all_models_both_operators_fewdata64_test.csv}
model, psnr-diff, ssim-diff, psnr-linblur, ssim-linblur
Small iResNet L=0.95,30.58464175996339,0.8926180309764351,24.756858298888282,0.7641911528338377
Small iResNet L=0.97,30.76602177006231,0.898464410969424,24.921227557615897,0.7750392112598838
Small iResNet L=0.99,30.890216415869613,0.9034084490912444,25.045522896382586,0.7831654390559151
Small iResNet L=0.995,30.974859217363214,0.9065141920532951,25.074223642980606,0.7852663604209785
Small iResNet L=0.999,30.94103602107676,0.9059878813561748,25.000289238150472,0.7828324971376597
iResNet L=0.95,30.479700579605744,0.8879496286289887,24.837590214813115,0.7697383621621595
iResNet L=0.97,30.816411358913165,0.8995971946069901,24.972573749633625,0.7791041491600146
iResNet L=0.99,30.946224590629082,0.9055967924836305,25.08697800259858,0.786034504644773
iResNet L=0.995,30.971466440231563,0.9059297028870643,25.062508226530525,0.7837725858079733
iResNet L=0.999,30.938613438940624,0.9056413964780312,25.036551410708604,0.7827335237441089
DiffNet,30.899892439610454,0.9069536700514921,25.011487367849256,0.784606321708649
U-Net,30.04179767591355,0.8861592514730687,24.082744175130824,0.7380826329184118
ConvResNet,29.75480257677614,0.8715528229693594,24.643252793120446,0.7548507996994018
\end{filecontents*}
\pgfplotstabletypeset[
    col sep=comma,
    columns/model/.style={string type, column name={ }, column type=l, 
    		},
    columns/psnr-diff/.style={column name={PSNR}, precision=2, column type/.add={|}{}},
    columns/ssim-diff/.style={column name={SSIM}, precision=3, column type/.add={}{}},
    columns/psnr-linblur/.style={column name={PSNR}, precision=2, column type/.add={|}{}},
    columns/ssim-linblur/.style={column name={SSIM}, precision=3, column type/.add={}{|}},
    assign column name/.style={/pgfplots/table/column name={\textbf{#1}}},
    every head row/.style={before row={ & \multicolumn{2}{c|}{\textbf{Nonlinear Diffusion}} & \multicolumn{2}{c|}{\textbf{Linear Blurring}}  \\}, after row=\hline},
	every row no 4/.style={after row=\hline},
	every row no 9/.style={after row=\hline},
	every row 8 column 1/.style={highlight bold},
	every row 3 column 1/.style={highlight bold},
	every row 3 column 2/.style={highlight bold},
	every row 10 column 2/.style={highlight bold},
	every row 7 column 3/.style={highlight bold},
	every row 7 column 4/.style={highlight bold},
	fixed zerofill=true
    ]{csv_files/all_models_both_operators_fewdata64_test.csv}
\end{table}


\subsection{Why invertibility matters}
\label{sec:investigate_reg_prop}

Unlike all other reconstruction methods considered in this paper, the iResNet is invertible, allowing us to analyze and interpret the behavior of the learned forward operator. Specifically, while $\varphi_\theta^{-1}$ has been trained for reconstruction (i.e., mapping $z^\delta$ to $x^\dagger$), we can simultaneously access $\varphi_\theta$ of the trained network (the inverse of the reconstruction scheme). This allows us to compare $\varphi_\theta$ with the true forward operator $F$ in several different ways.

\begin{figure}[t]
\centering
\includegraphics[width=0.95\textwidth]{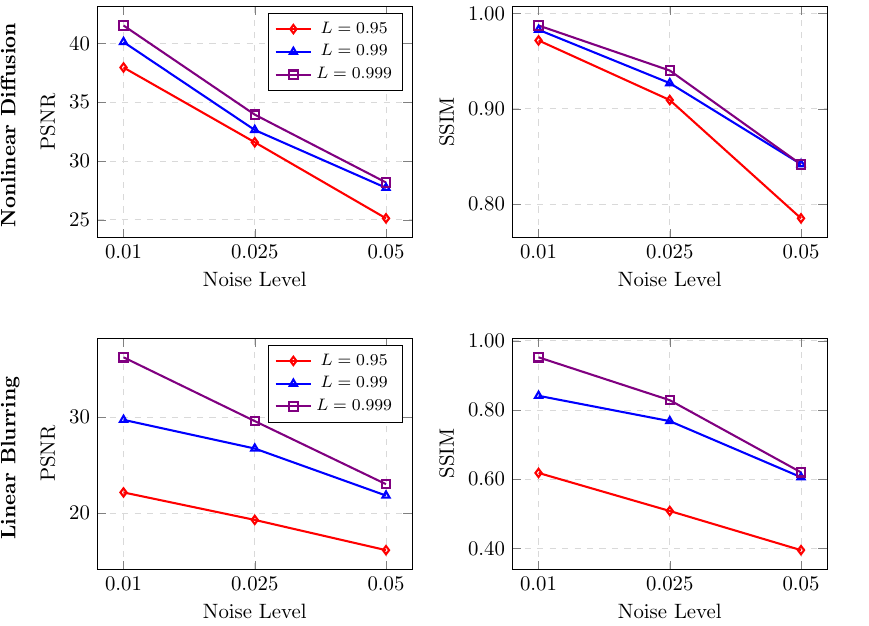}
\caption{Evaluation of the approximation accuracy of the true forward operator $F$ by iResNets $\varphi_\theta$, which were trained via reconstruction training (top row: \textbf{nonlinear diffusion}, bottom row: \textbf{linear blurring}). PSNR and SSIM values between $\varphi_\theta(x)$ and $F(x)$ are computed on the test dataset.
For the sake of clarity, we only illustrate three different values for the Lipschitz parameter $L$.}
\label{fig:approx_error}
\end{figure}

\begin{figure}[t]
\centering
\includegraphics[width=\textwidth]{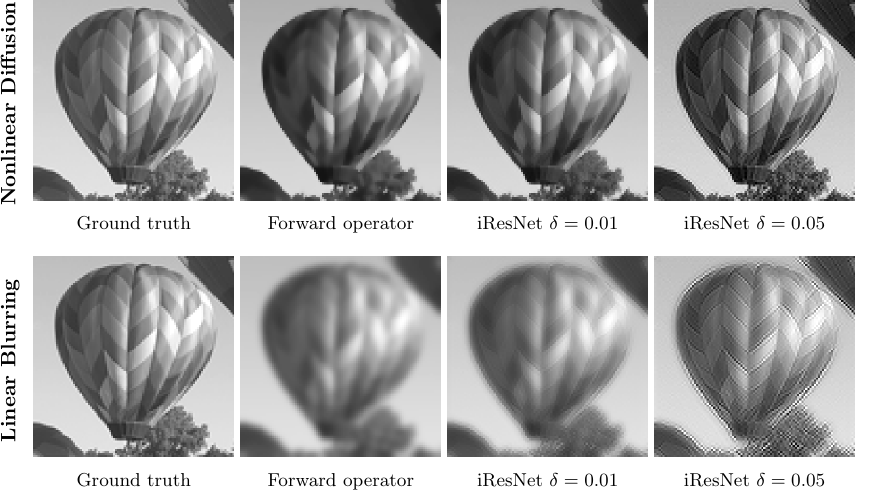}
\caption{Comparison of the forward operator (top row: \textbf{nonlinear diffusion}, bottom row: \textbf{linear blurring}) with the trained iResNets for different noise levels. From left to right: Ground truth image, blurred image using the forward operator with no added noise, as well as the output (forward pass) of the network trained with Lipschitz parameter $L=0.999$ for noise level $0.01$ and $0.05$. }\label{fig:forward_operator_net}
\end{figure}

To start, let us evaluate the accuracy of $\varphi_\theta$ in approximating $F$. Figure~\ref{fig:approx_error} illustrates the approximation performance -- measured in terms of PSNR and SSIM -- for both the nonlinear diffusion and linear blurring operators. The results indicate that as the noise level $\delta$ decreases, the approximation accuracy improves across all Lipschitz parameters $L$. Additionally, the highest Lipschitz parameter consistently yields the best approximation performance on average. This suggests that choosing a Lipschitz parameter close to one benefits both target reconstruction and forward operator approximation.

Examples of the true forward operator and the learned forward operator for $\delta=0.01$ as well as $\delta=0.05$ can be seen in Figure~\ref{fig:forward_operator_net}. One can observe that the learned forward operator resembles the true forward operator for small noise ($\delta=0.01$), displaying slightly less blurring or diffusion in both linear and nonlinear cases. In contrast, at a higher noise level ($\delta=0.05$), the learned forward operators tend to over-amplify details, especially edges in the images. The reason for this is that for high noise levels the learned operators need to exhibit a stronger regularization to be able to stably reconstruct the target images.

In the following sections, we explore this behavior in greater depth and draw connections to the learned regularization.
We would like to note that we present the results of our investigations solely for the nonlinear diffusion operator, as the results were comparable for the linear blurring operator.
Moreover, we emphasize that while this paper focuses on deblurring tasks, the proposed approaches are applicable to various problems and should be viewed as exemplary methods for numerically investigating the regularizing behavior of iResNets.

\subsubsection{Investigation of the learned forward operator by linearization}\label{sec:clustering}

To gain a deeper understanding of the learned forward operator $\phitheta$, we compare its local behavior to that of the true operator $F$. However, the inherent nonlinearity of both the network and the operator makes this analysis particularly challenging.

To address this challenge, we approximate $\phitheta$ and $F$ around a given image $x_0\in \R^{n\times m}$ with the help of a first order Taylor expansion. The first order Taylor expansion $T\Psi(\cdot,x_0): \R^{n\times m}  \to  \R^{n\times m}$ of a differentiable operator $\Psi:X\to X$ around $x_0$ restricted to some pixel $(n_0,m_0)\in \{ 1,\dots, n \}\times \{ 1,\dots, m \}$ is given by
\begin{equation*}
    T\Psi(x,x_0)\big\vert_{(n_0,m_0)} = \Psi(x_0)\big\vert_{(n_0,m_0)} + \left\langle J\Psi(x_0)\big\vert_{(n_0,m_0)}, x-x_0\right\rangle + R\Psi(x,x_0)\big\vert_{(n_0,m_0)}
\end{equation*}
with Jacobian 
\begin{equation*}
    J\Psi(x_0)\big\vert_{(n_0,m_0)} = \left(\frac{\partial {\Psi}_{n_0,m_0}}{\partial x_{k,l}} (x_0)\right)_{k=1,\dots,n,\, l=1,\dots,m},
\end{equation*}
$x\in \R^{n\times m}$ and remainder $R\Psi(\cdot,x_0): \R^{n\times m}  \to  \R^{n\times m}$.
In the setting of $\Psi$ being our iResNet or the nonlinear diffusion operator, the Jacobian can be interpreted as a linear blurring kernel that approximates the nonlinear diffusion operations of the network and the operator, respectively. This perspective allows us to assign a kernel to each pixel in the image $x_0$, providing a visual interpretation of $\phitheta$ and $F$. This approach was first introduced in \cite{simonyan_2014} in the context of image classification, leading to the creation of so-called saliency maps, which enable the interpretation of decisions made by image classification networks. While saliency maps in classification tasks are limited to the number of classes, the large number of image pixels in our setting makes manual analysis impractical. To address this, we opted to cluster the kernels of $\phitheta$ and $F$ and then analyze the resulting clusters, respectively.

Before discussing the clustering results, we will first provide an overview of our calculation pipeline.
The computation of the Jacobians and the clustering is performed separately for each image in the test set. This approach is necessary because we observed that clustering becomes ineffective or fails to converge when the number of pixels exceeds $\num{15000}$. We found that limiting the number of pixels per image to $\num{1500}$ leads to a good trade-off between cluster accuracy and computational efficiency.

We randomly select a set $\mathcal{P}$ of pixel indices, which remains the same across all images, networks, and the operator.
For each pixel $(n_0,m_0)$ in the set $\mathcal{P}$ and fixed image $x_0$, we calculate the Jacobians $J\phitheta(x_0)\big\vert_{(n_0,m_0)}$ and $JF(x_0)\big\vert_{(n_0,m_0)}$. For the network $\phitheta$, this can be achieved efficiently via backpropagation, giving access to the gradients of the input image. 
Each Jacobian is then cropped to a $9\times 9$ region centered around the pixel $(n_0,m_0)$. This cropping is performed because pixel values outside this region are typically near zero and do not contain relevant information. In what follows, the cropped image is referred to as the saliency map or linear blurring kernel. Before clustering, we normalize the values of these saliency maps while preserving the sign of the values.

Finally, the transformed saliency maps are clustered using spectral clustering. In addition to spectral clustering, we also evaluated other clustering techniques, ranging from primitive methods to more advanced approaches involving principal component analysis (PCA).
However, the clustering results were all very similar, which is why we opted for the standard and straightforward spectral clustering algorithm. 
This method is particularly effective when the clusters are expected to be highly non-convex. We refer the reader to \cite{luxburg_2007} for a more detailed explanation of spectral clustering. 

When applying a clustering method, selecting the number of clusters is a critical hyperparameter that must be determined in advance. Several techniques exist to help identify an appropriate number of clusters.
To obtain an initial estimate, we compare the predictions from the elbow method \cite{Thorndike_1953}, the gap statistic \cite{tibshirani_2001}, and the gap$^\ast$ statistic \cite{mohajer_2010}.
While the elbow method is heuristic in nature, both the gap statistic and the gap$^\ast$ statistic provide a more formal approach grounded in statistical analysis. In the examples presented in this section, we found that using two clusters offers a reasonable trade-off between achieving meaningful data separation and maintaining comparability between the network and the operator.

\begin{figure}[t]
\centering
\includegraphics[width=\textwidth]{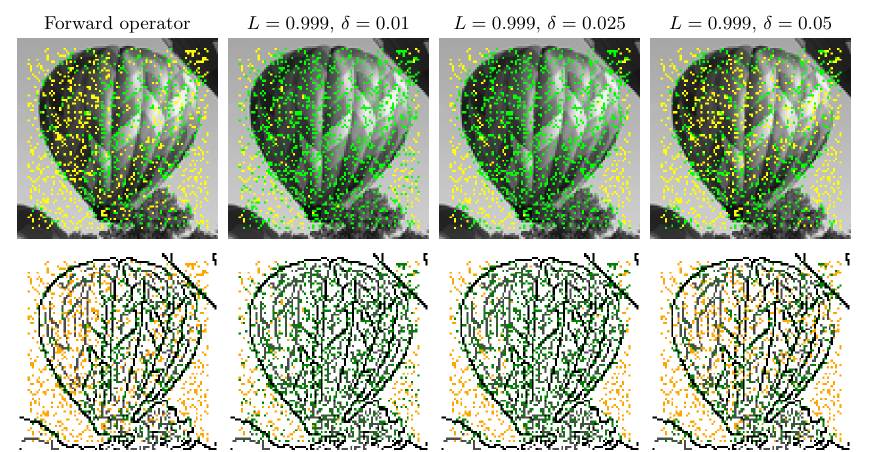}
\caption{Clustering of the saliency maps of the operator and the network trained with Lipschitz parameter $L=0.999$ for different noise levels. The first row visualizes the clustering with the ground truth image. The second row visualizes the clustering with the edges of the ground truth image. Weak edges are coloured in gray and strong edges are coloured in black.}\label{fig:clustering_different_nlvl}
\end{figure}

The clustering results for the operator as well as the network trained with different noise levels and fixed Lipschitz parameter $L=0.999$ are visualized in Figure~\ref{fig:clustering_different_nlvl}, exemplified on the balloon image from Figure~\ref{fig:forward_operator_net}. Each pixel for which saliency maps are calculated and clustered is colored based on its cluster assignment. We visualize the results alongside the ground truth image as well as its edges. In all subsequent plots, the cluster with the highest number of pixels that align with the image edges is highlighted in green.
The edges are calculated using the Canny edge detection algorithm. We differentiate between weak and strong edges, where strong edges are defined as those whose value exceeds $20\%$ of the maximum of the intermediate edge image generated during the Canny algorithm after applying Sobel filtering and gradient magnitude thresholding.

Figure~\ref{fig:clustering_different_nlvl} shows that the operator's saliency maps are primarily clustered into pixels corresponding to strong edges and those associated with weak edges and smooth areas. In contrast, the clustering for the network trained with $\delta=0.01$ indicates that the saliency maps are very similar for most pixels, as the majority of pixels fall into the cluster visualized in green. This behaviour is less pronounced for the network trained with $\delta=0.025$, where the saliency maps are more distinctly clustered into pixels corresponding to edges and those belonging to smoother regions of the image. The clustering of the network trained with $\delta=0.05$ most closely resembles that of the operator, with weak edges and smooth areas grouped into the same cluster. 

\begin{figure}[t]
\includegraphics[width=\textwidth]{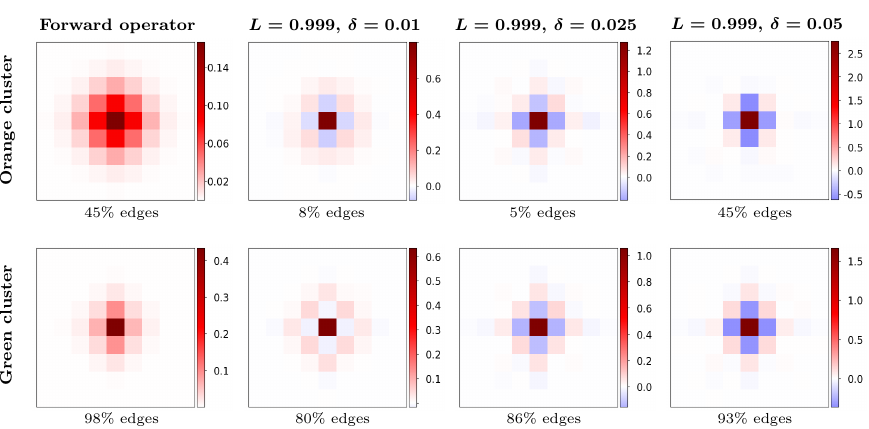}
\caption{Mean of the saliency maps of the operator and the network trained with Lipschitz parameter $L=0.999$ and noise levels $0.01$, $0.025$, $0.05$. The percentage of pixels on edges for each cluster is indicated below each saliency map.}\label{fig:saliency_maps_mean}
\end{figure}

To further investigate the clustering, we visualize the mean saliency maps for the individual clusters of the operator and the network trained with $L=0.999$ for different noise levels in Figure~\ref{fig:saliency_maps_mean}. 
The operator's mean saliency map of the orange cluster (comprising pixels on weak edges and in smooth areas) is much more spread out compared to that of the green cluster (strong edges). This was to be expected as the strength of the nonlinear diffusion depends on the gradient of the image, with stronger diffusion in areas with smaller gradients.

When comparing the network’s saliency maps to those of the operator, it is particularly striking that the spread of the orange cluster is significantly lower and decreases further as the noise level increases. Since the orange cluster primarily includes pixels in smooth regions and on weak edges, this leads to significantly less blurring in these areas and has a regularizing effect. Additionally, for small noise levels ($\delta=0.01, 0.025$), the dispersion of the network’s saliency maps is quite similar across both clusters, which may explain the poor clustering observed in Figure~\ref{fig:clustering_different_nlvl}. 

Furthermore, the network's mean saliency maps exhibit a strong emphasis on the central pixel, displaying significantly higher values compared to those of the operator. This is particularly the case for the cluster containing pixels in smooth areas as well as high noise levels. To counterbalance this and maintain the range of values in the blurred image, the values of the pixels adjacent to the central pixel are negative.

Overall, the behaviour of the saliency maps leads to an overemphasis of edges and significantly reduced blurring in smooth regions, especially for high noise levels. This aligns with the observations presented in Figure~\ref{fig:forward_operator_net}. 

We also examined the clustering results of the network with varying Lipschitz parameters. Our analysis confirms that smaller Lipschitz parameters significantly reduce the network's expressiveness, as previously discussed. The results can be found in Appendix~\ref{sec:cluster different lip}. 

So far, the investigations have been based entirely on a single image from the test set, which may lead to biased conclusions. This limitation arises because comparing clustering results across different images, let alone the entire test set, is challenging. Pixel cluster affiliations can vary significantly between different networks and images, making automated comparison difficult.
To obtain more reliable results, we decided to manually cluster the pixels of each image in the test set to create comparable clusters. Based on the operator's behavior and the clustering results of the balloon image, we chose to divide the pixels of each image into two groups: those belonging to smooth areas and those belonging to the edges. 

We implement this manual clustering using Canny edge detection to extract the edges of each image in the test set, thereby creating two distinct sets of pixels corresponding to smooth regions and edges.
To ensure a clear spatial separation between the two clusters, we dilate the edge image using a $3\times 3$ kernel, assigning pixels to the cluster representing smooth areas only if they are at least two pixels away from an edge. 
Moreover, the boundary region of the image is excluded from clustering to ensure that the extracted and cropped saliency maps consistently have the same size. For each image in the test set, we then sample $250$ pixels equidistantly from each cluster, calculate the corresponding saliency map, and average them within each cluster.

\begin{figure}[t]
\includegraphics[width=\textwidth]{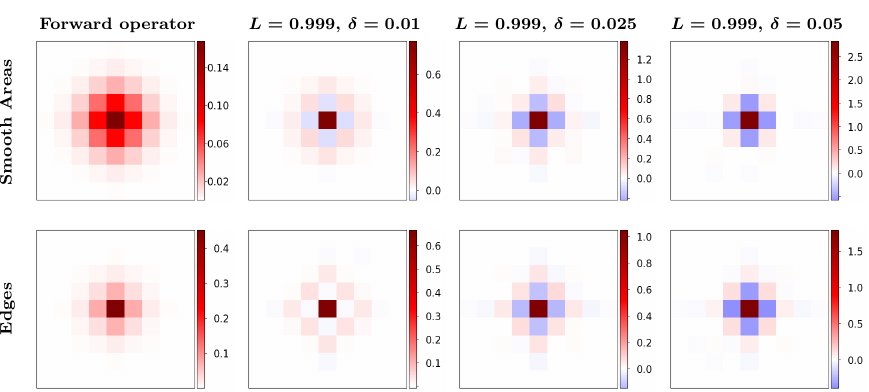}
\caption{Mean of the saliency maps of the operator and the network trained with Lipschitz parameter $L=0.999$ and noise levels $0.01$, $0.025$, $0.05$. The clustering is done manually on the test set according to edge pixels and smooth areas.}\label{fig:saliency_maps_mean_testset}
\end{figure}

The resulting averaged saliency maps are depicted in Figure~\ref{fig:saliency_maps_mean_testset} for the operator as well as the network trained with noise level $\delta=0.01, 0.025, 0.05$ and Lipschitz parameter $L=0.999$. The mean saliency maps exhibit a behavior remarkably similar to that observed in the clustering of the balloon image. This similarity confirms the reliability and significance of the balloon image's clustering results, indicating that this behavior generalizes well to other images.

\subsubsection{Investigation of regularization properties}\label{sec:directional_deriv}

\begin{figure}[t]
\centering
\includegraphics[width=0.94\textwidth]{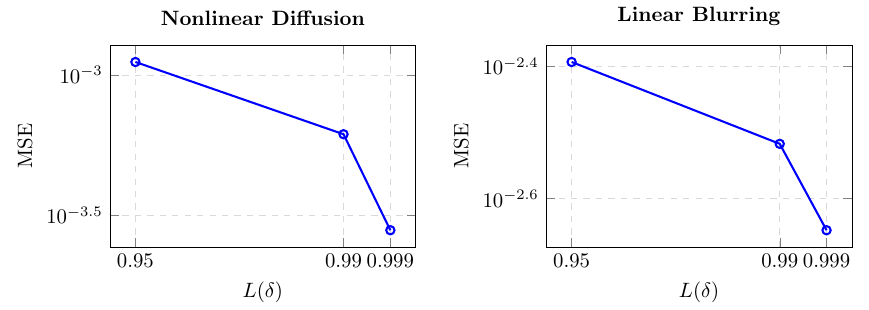}
\caption{Mean inversion error $\| \varphi_{\theta(L)}^{-1}(F(x)) - x \|$ on the test set for the nonlinear diffusion (left) and the linear blurring operator (right). The parameter choice $L(\delta)$ required for our convergence result in Lemma~\ref{lem:convergence} and Equation~\eqref{eq:inversion_error_convergence_prop} is given by $L(0.01)=0.999,L(0.025)=0.99,L(0.05)=0.95$. To be more precise, for $\delta=0.05$ the network was trained with $L=0.95$, for $\delta=0.025$ with $L=0.99$ and for $\delta=0.01$ with $L=0.999$. 
}\label{fig:inversion_error}
\end{figure}

As explained in Section~\ref{sec:theory}, iResNets ensure the existence and uniqueness of solutions $\varphi_\theta^{-1}(z^\delta)$, as well as the stability of the resulting reconstruction scheme with respect to $z^\delta$.
To finally qualify as a convergent regularization scheme, the reconstructions $\varphi_\theta^{-1}(z^\delta)$ must also converge to the ground truth $x^\dagger$ as the noise level $\delta$ approaches zero.
In Lemma~\ref{lem:convergence}, we derived a simple convergence result based on the inversion error $\| \varphi_{\theta(L)}^{-1}(F(x^\dagger)) - x^\dagger \|$ together with a suitable parameter choice rule $L(\delta)$ with $L(\delta)\to 1$ for $\delta\to 0$. To check whether the assumptions of the lemma are fulfilled, we evaluate the average inversion error on the test set, as illustrated in Figure~\ref{fig:inversion_error}.
It can be observed that the error decreases as $L(\delta) \to 1$ for both forward operators, indicating that the iResNet reconstruction approach meets the conditions for a convergent regularization scheme on average.

\begin{figure}[p]
\includegraphics[width=\textwidth]{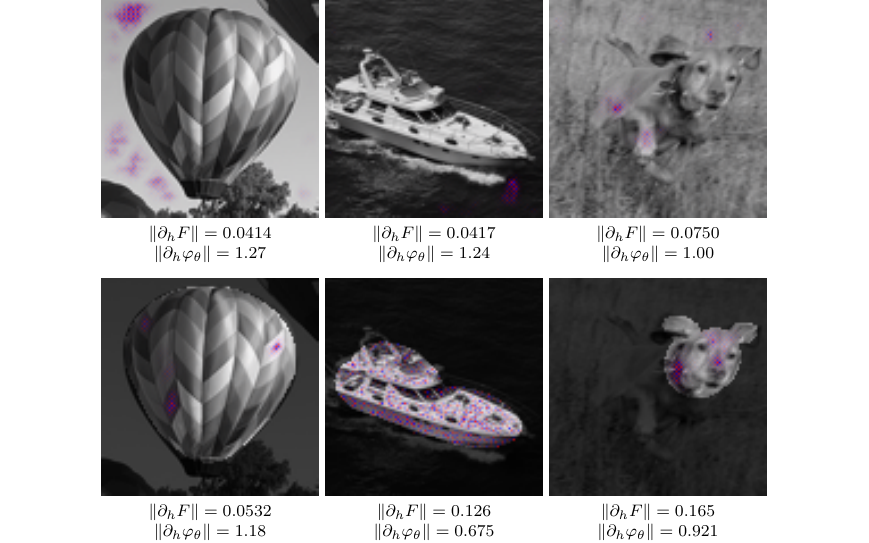}
\caption{Computation of directions $h$ for which the difference $\| \partial_h \varphi_\theta(x_0) \| - \| \partial_h F(x_0) \|$ is particularly large ($\varphi_\theta$ trained for the \textbf{nonlinear diffusion operator} with $\delta=0.01$). The direction vector $h$ is displayed in color on the image $x_0$ (grayscale). It can be interpreted as a direction in which the network has learned a significant regularization.
In the bottom row, $h$ is restricted to the subject of the image $x_0$.}
\label{fig:reg_direction}

\centering

\includegraphics{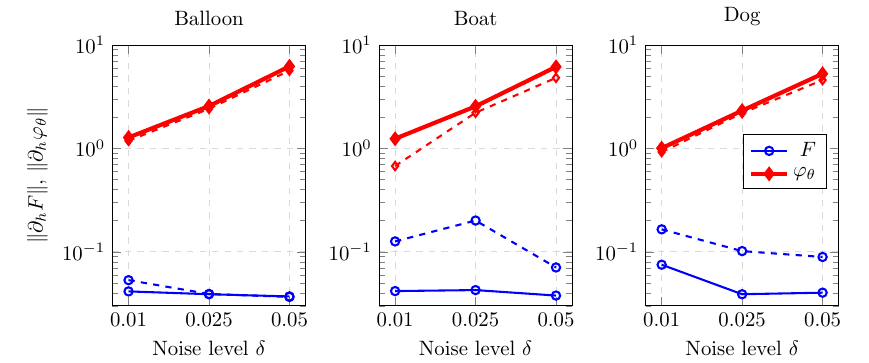}
\caption{Computation of directions $h$ for which the difference $\| \partial_h \varphi_\theta(x_0) \| - \| \partial_h F(x_0) \|$ is particularly large ($\varphi_\theta$ trained for the \textbf{nonlinear diffusion operator} with different noise levels $\delta$). For the three test images ``Balloon'', ``Boat'' and ``Dog'', we display $\| \partial_h \varphi_\theta(x_0) \|$ in red and $\| \partial_h F(x_0) \|$ in blue. The dashed lines show the results, where $h$ is restricted to the subject of the image (cf.\ Figure \ref{fig:reg_direction}).}
\label{fig:reg_dir_plot}
\end{figure}

In addition to the general regularization properties, the invertibility of iResNets allows for a deeper analysis of how $\varphi_\theta$ regularizes the forward operator $F$. For this purpose, we consider the concept of local ill-posedness and compare $\varphi_\theta$ with $F$.

Recall that a nonlinear inverse problem is called locally ill-posed in a point $x_0 \in X$, if any open neighborhood of $x_0$ contains a sequence $(x_k)$ such that $F(x_k) \to F(x_0)$ but $x_k \not\to x_0$ \cite[Definition~3.15]{Schuster_2012}. 
For differentiable $F$, this often entails that some directional derivatives $\partial_h F(x_0) = F'(x_0)h$ are small. Accordingly, for our nonlinear diffusion operator $F$ we expect $\| \partial_h F(x_0) \|$ to be small for certain directions $h \in X$, $\|h\|=1$, dependent on the input image $x_0 \in X$.

In contrast, we expect the trained iResNet $\varphi_\theta$ to have learned a regularization for the inverse problem. Thus, we anticipate $\| \partial_h \varphi_\theta(x_0)\| \gg \| \partial_h F(x_0)\|$ for the same directions $h$. The Lipschitz constraint even ensures that $\| \partial_h \varphi_\theta(x_0)\| \geq 1-L$ is guaranteed.
Moreover, we can use gradient ascent to find directions $h \in X$, $\|h\|=1$ such that
\begin{equation*}
    \|\partial_h \varphi_\theta(x_0)\|^2- \|\partial_h F(x_0)\|^2
\end{equation*}
is as large as possible. In other words, we look for a direction $h$ in which the network has learned the strongest regularization with respect to a given input image $x_0$. We can also restrict $h$ to certain areas of the image $x_0$ in order to analyze the regularization inside this area.

Figure \ref{fig:reg_direction} shows the resulting directions $h$ of highest regularization for three different test images $x_0$. Note that gradient ascent does not lead to unique solutions, and the directions should therefore rather be seen as examples.

The first striking observation is that $h$ consistently exhibits a checkerboard pattern. This occurs because adding a small perturbation $\varepsilon h$ to $x_0$ introduces numerous small edges, which are subsequently blurred by the forward operator $F$. Since these small edges may represent important image details, regularization by the iResNet is necessary. Additionally, $h$ is primarily concentrated in the smooth parts of the images, such as the uniform areas of the background. In these regions, the introduction of small edges has the most pronounced effect. Moreover, reconstructing a smooth background without any details is a relatively simple task for a neural network, requiring only a few information from the noisy and blurred image $z^\delta$ and relying more heavily on the learned prior. Therefore, it is optimal to apply strong regularization in these areas. Even when $h$ is restricted to the subject of the images, it remains concentrated in the smoother regions.

An iResNet trained on higher noise levels is expected to apply a stronger regularization. To verify this, we compare the values of $\|\partial_h \varphi_\theta(x_0)\|$ for $\varphi_\theta$ trained on different noise levels $\delta$. The result is illustrated in Figure \ref{fig:reg_dir_plot}. For all three test images, we indeed observe that the amount of regularization increases with $\delta$.

 
\section{Discussion and Outlook}\label{sec:outlook}

This present work builds on the theoretical foundation of the iResNet reconstruction approach for solving inverse problems, previously introduced in \cite{iresnet_01_regtheory, iresnet_02_bayesian}, by conducting extensive numerical experiments on two real-world tasks. Our experiments demonstrate the competitiveness of iResNets compared to state-of-the-art reconstruction methods. Earlier works \cite{iresnet_01_regtheory, iresnet_02_bayesian} established criteria for rendering the iResNet reconstruction approach a convergent regularization scheme. However, their efficacy in real-world, high-dimensional tasks remained unverified, and initial numerical experiments revealed that the proposed architectures did not scale well to such problems. To address this, we developed a large-scale iResNet architecture and evaluated its performance on both a linear blurring problem and a nonlinear diffusion problem.

Our results show that iResNets achieve performance comparable to leading image reconstruction networks, such as U-Nets. However, unlike most state-of-the-art methods, iResNets ensure the stability of the resulting reconstruction scheme by controlling the hyperparameter $L$. This leads to increased robustness against adversarial attacks, varying noise during testing, and improved performance in scenarios with limited training data -- surpassing classical network architectures in these aspects.
Additionally, the inherent invertibility of iResNets allows for deeper insights into the learned forward operator and the characteristics of the learned regularization. This feature has the potential to reduce the black-box nature of learned reconstruction methods, enhancing interpretability and trustworthiness. Despite these benefits, iResNets come with a trade-off: significantly longer training times compared to state-of-the-art methods.

In our numerical experiments, the best trade-off between interpretability of the learned solution and performance was achieved by DiffNets, as introduced in \cite{arridge_hauptmann_2020}, though these are specifically tailored for diffusion problems. In contrast, our iResNet reconstruction approach is applicable to a broader class of inverse problems and offers theoretical guarantees in addition to interpretability. On this note, it would be valuable to explore whether similar techniques as in the construction of DiffNets could be applied to iResNets. In particular, designing the architecture to fit the specific problem could help reduce the number of network parameters and speed up training.
More generally, reducing iResNet training times while maintaining performance is crucial for improving their practical applicability. Possible directions include allowing for larger gradient steps -- currently constrained by the Lipschitz condition -- and refining the computation of the Lipschitz constant during training.

To fully assess the versatility of iResNets, future research should investigate their performance across a wider range of forward operators. Additionally, exploring the effectiveness of iResNets as denoisers in Plug-and-Play (PnP) schemes could pave the way for a new range of applications for these models.

In summary, this work provides initial numerical evidence supporting the iResNet regularization scheme as an effective and interpretable learned reconstruction method with theoretical guarantees. Consequently, iResNets offer a promising approach to solving complex real-world inverse problems by bridging the gap between high-performance, data-driven methods that often lack theoretical justification and classical methods with guarantees but inferior performance. This sets the stage for further research in this area.

\newpage
\appendix

\section{Transforming the inverse of iResNets into averaged operators}
\label{sec:iResNet_averaged}

In the context of Plug-and-Play (PnP) algorithms, averaged operators play a central role in establishing convergence results. Recall that an operator $A$ is $t$-averaged if it can be expressed as $A=(1-t)\,\text{Id} + t R$ for some operator $R$ satisfying $\text{Lip}(R)\leq 1$ and some $t\in (0,1)$.
The introduced iResNet, and specifically its inverse, generally does not satisfy this property. However, one can define an averaged surrogate operator by leveraging a fixed oracle element $\hat{x}\in X$. This oracle could, for instance, be a denoised version of the noisy input image, generated using a classical and computationally inexpensive denoiser.
As this work uses the inverse  $\varphi_\theta^{-1}$ to address the inverse problem, our analysis will focus on $\varphi_\theta^{-1}$ in the following discussion.
For simplicity, we assume that $\varphi_\theta$ consists of a single subnetwork, i.e., $\varphi_\theta = \text{Id} - f_\theta$. The inverse can then be expressed as $\varphi_\theta^{-1} = \nicefrac{1}{1-L^2}\, (\text{Id} - g)$ for some $g: X\to X$ with $\text{Lip}(g)\leq L < 1$ \cite[Remark 2.3]{iresnet_02_bayesian}. 
Now, consider the operator
\begin{equation*}
    D(x) = L \hat{x} + (1-L)\,\varphi_\theta^{-1}(x) =  L\, \hat{x} + \frac{1}{1+L} (x-g(x)) \quad \text{for } x\in X.
\end{equation*}
Following the same line of reasoning as in \cite[Lemma 18]{Hertrich2023PnP}, one can show that $D$ is $\nicefrac{1}{2}$-averaged. The proof relies on the observation that $g$ is $t$-averaged with $t = \nicefrac{L+1}{2}$, which follows directly from the definition of averagedness and the Lipschitz continuity of $g$, cf.\ \cite[Theorem 2]{Hertrich2023PnP}. 

\section{Proof of Lemma~\ref{lem:convergence}}\label{Sec:Proof_convergence}

\begin{proof} 
The Lipschitz continuity of the inverse, cf.\ Equation \ref{Eq:Lip concat inverse}, directly implies
\begin{align*}
     \| \varphi_{\theta(L(\delta))}^{-1}(z^\delta) - x^\dagger\| &\leq  \| \varphi_{\theta(L(\delta))}^{-1}(z^\delta) - \varphi_{\theta(L(\delta))}^{-1}(F(x^\dagger)) \| + \| \varphi_{\theta(L(\delta))}^{-1}(F(x^\dagger)) - x^\dagger \| \\
     &\leq \frac{1}{1-L(\delta)} \|z^\delta - F(x^\dagger)\| + \| \varphi_{\theta(L(\delta))}^{-1}(F(x^\dagger)) - x^\dagger \| .
\end{align*}
The last inequality yields the desired assertion due to \eqref{eq:inversion_error_convergence_prop} and \eqref{eq:conv_param_choice_ass}.
\end{proof}

\section{Proof of Lemma \ref{lem:Sh_fne}}
\label{sec:proof_Sh_fne}

\begin{proof}
We denote $u_{t+1} = S_{g_\lambda,h}(u_t)$ and $v_{t+1} = S_{g_\lambda,h}(v_t)$ to enhance readability of the notation. It holds
\begin{equation*}
\begin{split}
\|(u_t - u_{t+1}) &- (v_t - v_{t+1})\|^2 = \|(u_t - v_t) - (u_{t+1} - v_{t+1})\|^2 \\
&= \|u_t - v_t\|^2 - 2 \langle u_t - v_t, u_{t+1} - v_{t+1}\rangle + \|u_{t+1} - v_{t+1}\|^2.
\end{split}
\end{equation*}
By \eqref{eq:pm_implicit_euler}, it follows
\begin{equation*}
\begin{split}
&\quad \,\, u_t - v_t \\
&= \left(u_{t+1} - h \cdot \mathrm{div} \left( g_\lambda(| \nabla u_{t+1} |) \nabla u_{t+1} \right) \right) - \left( v_{t+1} - h \cdot \mathrm{div} \left( g_\lambda(| \nabla v_{t+1} |) \nabla v_{t+1} \right) \right) \\
&= \left( u_{t+1} - v_{t+1} \right) - h \cdot \mathrm{div}\left(  g_\lambda(| \nabla u_{t+1} |) \nabla u_{t+1} - g_\lambda(| \nabla v_{t+1} |) \nabla v_{t+1}\right).
\end{split} 
\end{equation*}
Using this and the boundary condition (zero Dirichlet or zero Neumann), we obtain
\begin{equation*}
\begin{split}
&\quad \,\, -2 \langle u_t - v_t, u_{t+1} - v_{t+1}\rangle \\
 &= - 2 \| u_{t+1} - v_{t+1}\|^2 + 2 h \langle \mathrm{div}\left(  g_\lambda(| \nabla u_{t+1} |) \nabla u_{t+1} - g_\lambda(| \nabla v_{t+1} |) \nabla v_{t+1}\right), u_{t+1} - v_{t+1}\rangle\\
 &= - 2 \| u_{t+1} - v_{t+1}\|^2  - 2 h \langle  g_\lambda(| \nabla u_{t+1} |) \nabla u_{t+1} - g_\lambda(| \nabla v_{t+1} |) \nabla v_{t+1} , \nabla u_{t+1} - \nabla v_{t+1} \rangle.
\end{split}
\end{equation*}
This now implies
\begin{equation*}
\begin{split}
\|(u_t - u_{t+1}) &- (v_t - v_{t+1})\|^2 = \|u_t - v_t\|^2 - \|u_{t+1} - v_{t+1}\|^2 \\
&- 2 h \langle  g_\lambda(| \nabla u_{t+1} |) \nabla u_{t+1} - g_\lambda(| \nabla v_{t+1} |) \nabla v_{t+1} , \nabla u_{t+1} - \nabla v_{t+1} \rangle.
\end{split}
\end{equation*}
By Lemma~\ref{lem:pm_g_monotone}, the scalar product in the last term is greater than or equal to zero due to $|\nabla u_t|, |\nabla v_t| \leq \lambda$. Therefore, the proof is complete.
\end{proof}

\begin{lemma}
\label{lem:pm_g_monotone}
For some $\lambda > 0$, let $g_\lambda \in C^1(\R_{\geq 0})$ be a Perona-Malik filter function, which fulfills \eqref{eq:pm_filter_monotone}, and let $v, \tilde{v} \in L^2(\Omega, \R^n)$ be two functions with $|v|, |\tilde{v}| \leq \lambda$ almost everywhere. Then, it holds
\begin{equation*}
\langle g_\lambda(|v|) v - g_\lambda(|\tilde{v}|) \tilde{v}, v - \tilde{v} \rangle \geq 0. 
\end{equation*}
\end{lemma}
\begin{proof}
By definition, it holds
\begin{equation*}
\begin{split}
\langle g_\lambda(|v|) v - g_\lambda(|\tilde{v}|) \tilde{v}, v - \tilde{v} \rangle
&= \int_{\Omega}  \langle g_\lambda(|v(x)|) v(x) - g_\lambda(|\tilde{v}(x)|) \tilde{v}(x), v(x) - \tilde{v}(x) \rangle \, \mathrm{d}x.
\end{split}
\end{equation*}
We now show that the function under the integral is greater than or equal to zero. For arbitrary $w, \tilde{w} \in \R^n$ with $|w|, |\tilde{w}| \leq \lambda$, it holds
\begin{equation*}
\begin{split}
&\quad \,\, \langle g_\lambda(|w|) w - g_\lambda(|\tilde{w}|) \tilde{w}, w - \tilde{w} \rangle\\
 &=  \langle g_\lambda(|w|) w, w \rangle - \langle g_\lambda(|w|) w , \tilde{w} \rangle - \langle g_\lambda(|\tilde{w}|) \tilde{w}, w \rangle + \langle g_\lambda(|\tilde{w}|) \tilde{w}, \tilde{w} \rangle \\
&\geq  g_\lambda(|w|) |w|^2 - g_\lambda(|w|) |w| |\tilde{w}| -  g_\lambda(|\tilde{w}|) |w| |\tilde{w}| + g_\lambda(|\tilde{w}|) |\tilde{w}|^2 \\
&= \left( g_\lambda(|w|)|w| - g_\lambda(|\tilde{w}|)|\tilde{w}| \right) \left(|w| - |\tilde{w}|\right).
\end{split}
\end{equation*}
The last term is non-negative due to the monotonicity of $s \mapsto g_\lambda (s)s$.
\end{proof}


\section{Out-of-distribution data}\label{sec:ood_data}

To evaluate the performance of the reconstruction methods on out-of-distribution data, we test the networks and TV, trained on the STL-10 dataset, using the Low-Dose Parallel Beam (LoDoPaB)-CT dataset \cite{Leuschner2021Lodopab}. This dataset comprises human chest computed tomography (CT) scans, from which we utilize 128 ground-truth reconstructions from the test set. To align the CT reconstructions with the scale and dimensions of STL-10 images, we rescale and resize them to match the settings used in this study.

The reconstruction performance for the nonlinear diffusion and the linear blurring operator, with a noise level of $\delta=0.025$, is shown in Table~\ref{tab:performance_ood} for networks trained on $\num{16384}$ STL-10 images, and in Table~\ref{tab:performance_ood_fewdata} for networks trained on $\num{64}$ STL-10 images (see Section~\ref{sec:small_dataset}).

For the larger training dataset, the reconstruction methods exhibit performance comparable to that observed on the STL-10 test dataset, with the ConvResNet consistently outperforming all other methods across both forward operators. This behavior is likely due to the high diversity of images in the STL-10 dataset, making the LoDoPaB-CT images appear less out-of-distribution.

\begin{figure}[p]
\centering
\includegraphics[width=\textwidth]{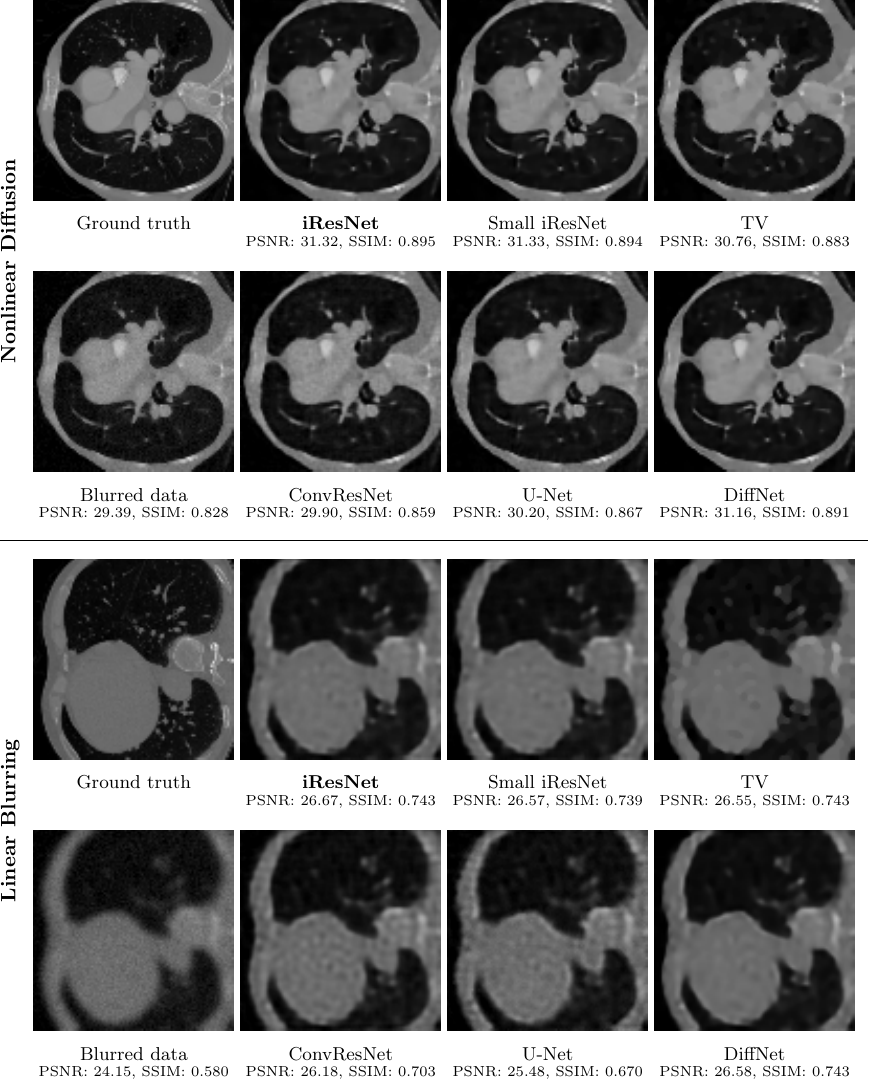}

\caption{Comparison of reconstructions for both forward operators with $\delta = 0.025$ on the LoDoPaB-CT dataset. All models were optimized on 64 samples from the STL-10 dataset with $\delta = 0.025$.}\label{fig:ood_lodopab}
\end{figure}

In contrast, for the smaller training dataset of $64$ samples, the iResNets significantly outperform both the U-Net and the ConvResNet, followed by the DiffNet and TV. These results suggest that in scenarios with limited training data, iResNets demonstrate greater robustness to outliers in the testing data. Figure~\ref{fig:ood_lodopab} provides examples of reconstructions for both forward operators using networks trained on 64 samples.

\begin{table}[h]
\centering
    \caption{Reconstruction performance for both forward operators with noise level $\delta = 0.025$ on 128 samples from the LoDoPaB-CT test set. The models were optimized on the STL-10 dataset with $\delta = 0.025$.}
    \label{tab:performance_ood}
    \begin{filecontents*}{csv_files/all_models_both_operators_ood_full_dataset_nlvl_0.025.csv}
model, psnr-diff, ssim-diff, psnr-linblur, ssim-linblur
Small iResNet L=0.95, 31.621425435281164, 0.8266537973153503, 27.05220147205776, 0.7099884176106344
Small iResNet L=0.97, 31.766623535211227, 0.8311864618233398, 27.235187415286294, 0.7205010353716914
Small iResNet L=0.99, 31.873003092325302, 0.8338163114854934, 27.61798907445445, 0.7379903886208428
Small iResNet L=0.995, 31.923337661085025, 0.835364381317653, 27.699504205160142, 0.740453186415493
Small iResNet L=0.999, 32.02111435406147, 0.8372109569002082, 27.771058637596045, 0.7425722923857487
iResNet L=0.95, 31.514293100206693, 0.8223626871104209, 27.09833994570263, 0.711727827528382
iResNet L=0.97, 31.744242716403274, 0.8304228114227356, 27.306571126702764, 0.7233056979010563
iResNet L=0.99, 31.916749319474388, 0.8347518109584682, 27.764421705607965, 0.7420490202080312
iResNet L=0.995, 32.01815565080331, 0.8375369858850559, 27.90160064821646, 0.7454787332472436
iResNet L=0.999, 32.121122607072174, 0.8387909040876009, 28.022736094365023, 0.7491601857318015
DiffNet, 32.139615826494556, 0.8387057338781326, 27.90796629183762, 0.7448871178563305
U-Net, 32.087603811005636, 0.8429722084511674, 26.762431394067853, 0.7320284454889457
ConvResNet, 32.35755272465627, 0.8465821022285942, 28.226015375814903, 0.7558094322892224
TV, 31.320772097626264, 0.8227659294129432, 27.25212490260622, 0.726024165236208
\end{filecontents*}
\pgfplotstabletypeset[
    col sep=comma,
    columns/model/.style={string type, column name={ }, column type=l, 
    		},
    columns/psnr-diff/.style={column name={PSNR}, precision=2, column type/.add={|}{}},
    columns/ssim-diff/.style={column name={SSIM}, precision=3, column type/.add={}{}},
    columns/psnr-linblur/.style={column name={PSNR}, precision=2, column type/.add={|}{}},
    columns/ssim-linblur/.style={column name={SSIM}, precision=3, column type/.add={}{|}},
    assign column name/.style={/pgfplots/table/column name={\textbf{#1}}},
    every head row/.style={before row={ & \multicolumn{2}{c|}{\textbf{Nonlinear Diffusion}} & \multicolumn{2}{c|}{\textbf{Linear Blurring}}  \\}, after row=\hline},
	every row no 4/.style={after row=\hline},
	every row no 9/.style={after row=\hline},
	every row 12 column 1/.style={highlight bold},
	every row 12 column 2/.style={highlight bold},
	every row 12 column 3/.style={highlight bold},
	every row 12 column 4/.style={highlight bold},
	fixed zerofill=true
    ]{csv_files/all_models_both_operators_ood_full_dataset_nlvl_0.025.csv}
\end{table}
\begin{table}[h] 
\centering
    \caption{Reconstruction performance for both forward operators with noise level $\delta = 0.025$ on 128 samples from the LoDoPaB-CT test set. The models were optimized on 64 samples from the STL-10 dataset with $\delta = 0.025$. }
    \label{tab:performance_ood_fewdata}
    \begin{filecontents*}{csv_files/all_models_both_operators_ood_fewdata64_nlvl_0.025.csv}
model, psnr-diff, ssim-diff, psnr-linblur, ssim-linblur
Small iResNet L=0.95, 31.169162442508004, 0.8112410454840806, 26.816254003057495, 0.6954755850342738
Small iResNet L=0.97, 31.392618773517356, 0.8197588527653896, 27.073629458957306, 0.7120654042841091
Small iResNet L=0.99, 31.54539666025093, 0.8264888439315693, 27.262162135153737, 0.7242669399211467
Small iResNet L=0.995, 31.649270856690663, 0.8296071930826975, 27.27643518272726, 0.7269409413149079
Small iResNet L=0.999, 31.66293967374921, 0.8301403307257202, 27.241276174381806, 0.7243359984586504
iResNet L=0.95, 31.06098910318032, 0.8063078887942161, 26.931039145364537, 0.703272134731297
iResNet L=0.97, 31.454017624129225, 0.8214170639296234, 27.110066716162855, 0.7156710017056462
iResNet L=0.99, 31.678243647663496, 0.8302858690909156, 27.306739939585462, 0.7278275669950875
iResNet L=0.995, 31.692142035616094, 0.8309461870898182, 27.28528203894109, 0.726673302854411
iResNet L=0.999, 31.689500065381118, 0.831071548055136, 27.274449593900215, 0.7250006119173658
DiffNet, 31.585270564372458, 0.8264461707189121, 27.203097022277635, 0.7246725718455114
U-Net, 30.314812007688204, 0.7900420634158484, 24.98124128508816, 0.642063229045498
ConvResNet, 30.167269234423117, 0.7812206429340963, 26.72217052468536, 0.6894777857909461
TV, 31.320772097626264, 0.8227659294129432, 27.252125525269975, 0.7260241530142392
\end{filecontents*}
\pgfplotstabletypeset[
    col sep=comma,
    columns/model/.style={string type, column name={ }, column type=l, 
    		},
    columns/psnr-diff/.style={column name={PSNR}, precision=2, column type/.add={|}{}},
    columns/ssim-diff/.style={column name={SSIM}, precision=3, column type/.add={}{}},
    columns/psnr-linblur/.style={column name={PSNR}, precision=2, column type/.add={|}{}},
    columns/ssim-linblur/.style={column name={SSIM}, precision=3, column type/.add={}{|}},
    assign column name/.style={/pgfplots/table/column name={\textbf{#1}}},
    every head row/.style={before row={ & \multicolumn{2}{c|}{\textbf{Nonlinear Diffusion}} & \multicolumn{2}{c|}{\textbf{Linear Blurring}}  \\}, after row=\hline},
	every row no 4/.style={after row=\hline},
	every row no 9/.style={after row=\hline},
	every row 8 column 1/.style={highlight bold},
    every row 9 column 1/.style={highlight bold},
	every row 8 column 2/.style={highlight bold},
    every row 9 column 2/.style={highlight bold},
	every row 7 column 3/.style={highlight bold},
	every row 7 column 4/.style={highlight bold},
	fixed zerofill=true
    ]{csv_files/all_models_both_operators_ood_fewdata64_nlvl_0.025.csv}
\end{table}

\clearpage

\clearpage
\section{Test-time noise level variations (additional material)}\label{App:varying_nlvl}%
\begin{figure}[H]
\includegraphics[width=0.95\textwidth]{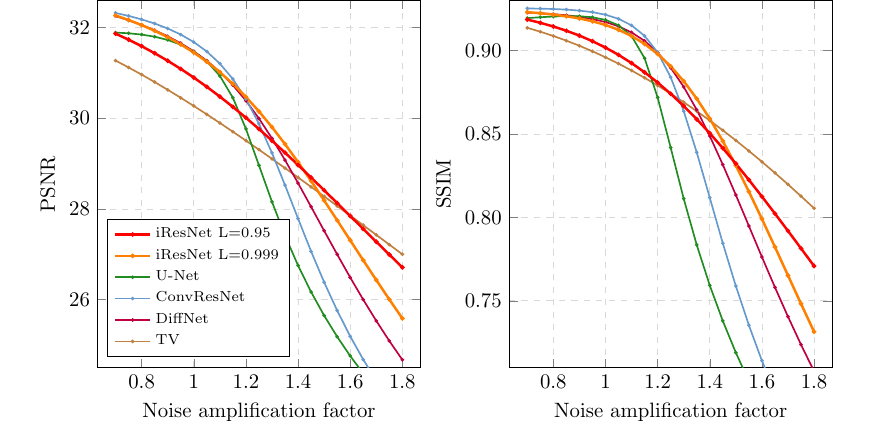}
\caption{Reconstruction performance on test data with noise level $\delta_{\text{test}} = c \cdot 0.025$ for the \textbf{nonlinear diffusion operator} and noise amplification factor $c \in [0.7, 1.8]$. All networks and hyperparameters were optimized on training data with $\delta=0.025$.}
\end{figure}

\begin{figure}[H]
\includegraphics[width=0.95\textwidth]{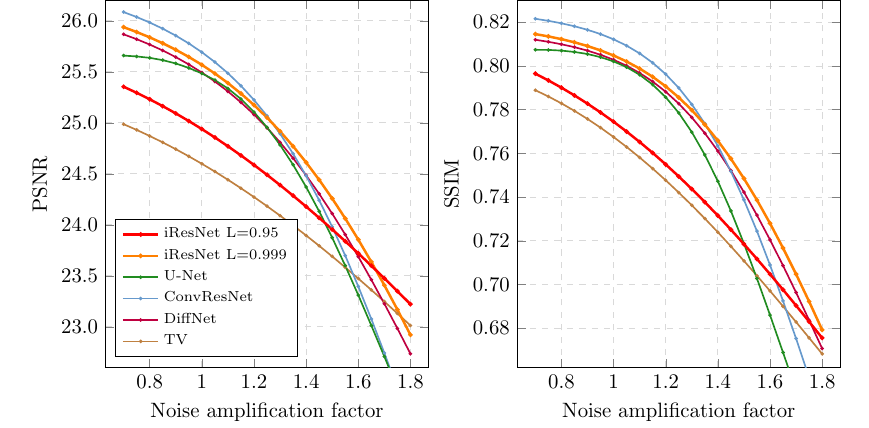}
\caption{Reconstruction performance on test data with noise level $\delta_{\text{test}} = c \cdot 0.025$ for the \textbf{linear blurring operator} and noise amplification factor $c \in [0.7, 1.8]$. All networks and hyperparameters were optimized on training data with $\delta=0.025$.}
\end{figure}


\newpage
\section{Cluster comparison for different Lipschitz parameters}\label{sec:cluster different lip}
\begin{figure}[h!]
\includegraphics[width=\textwidth]{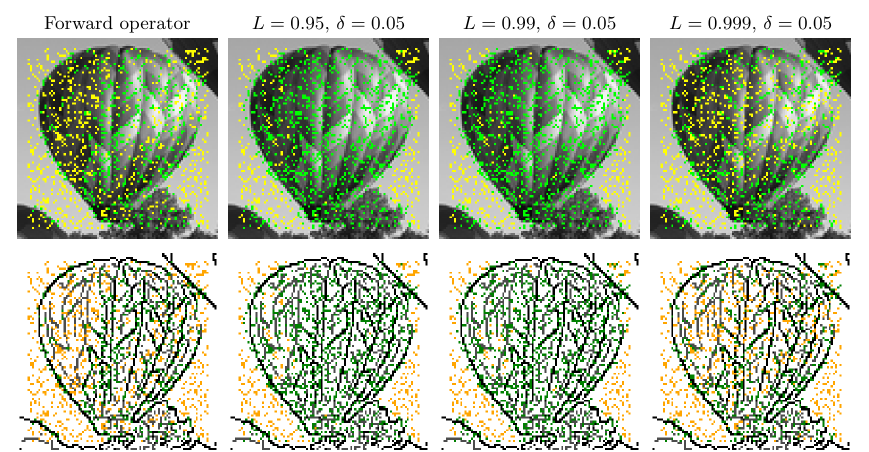}
\caption{Clustering of the saliency maps of the operator and the network trained with noise level $\delta=0.05$ for different Lipschitz parameters. The first row visualizes the clustering with the ground truth image and the second row with the corresponding edges (weak edges are gray and strong edges are black). }\label{fig:clustering_different_lip}

\includegraphics[width=\textwidth]{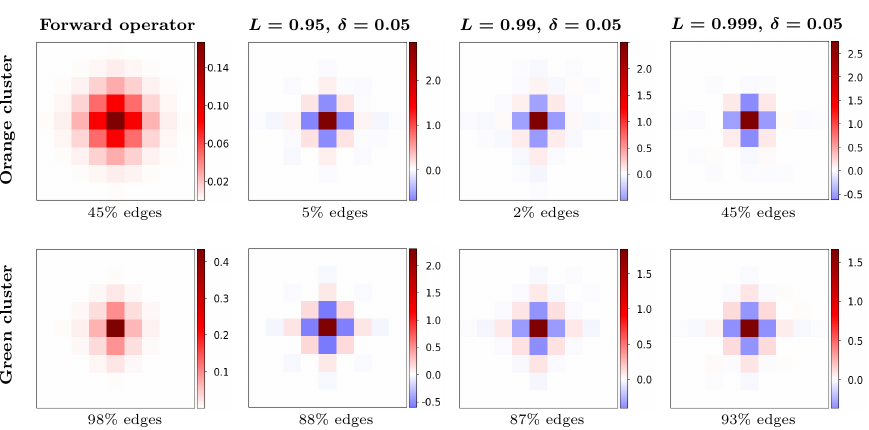}
\caption{Mean of the saliency maps of the operator and the network trained with Lipschitz parameters $L=0.95,0.99,0.999$ and noise level $0.05$. The percentage of pixels on edges for each cluster is indicated below each saliency map.}\label{fig:saliency_maps_mean_lipschitz}
\end{figure}

\bibliographystyle{siamplain}
\bibliography{literature}

\end{document}